\documentclass[11pt]{article}

\usepackage{a4wide}
\usepackage{amsmath,amsfonts,amssymb}

\usepackage{longtable}
\usepackage{comment}
\usepackage{manfnt}
\usepackage{xcolor}
\usepackage[all,cmtip]{xy}
\usepackage{mathtools}

\def\qed{\hfill $\Box$}
\def\gr{\rm gr}
\def\g{\gamma}

\def\pf{{\em Proof.\ }}

\def\<{\langle}
\def\>{\rangle}
\def\Z{\mathbb Z}
\def\k{\kappa}
\def\l{\lambda}
\def\n{\{1,\ldots, n\}}
\def\l{\ell}

\def\pf{{\noindent \it Proof. \ }}
\def\qed{\hfill $\Box$}
\def\g{\gamma}
\def\G{\Gamma}

\def\t{\overline{t}}

\def\<{\langle}
\def\>{\rangle}

\def\Z{\mathbb Z}    

\def\l{\lambda}

\newtheorem{theorem}{Theorem}[section]
\newtheorem{corollary}[theorem]{Corollary}
\newtheorem{lemma}[theorem]{Lemma}

\newtheorem{proposition}[theorem]{Proposition}

\title{The Formanek-Procesi group with base a right-angled Artin group: Lie algebra and Residual nilpotence}
\author{V. Metaftsis and A.I. Papistas}
\date{}

\begin{document}

\maketitle

\begin{abstract}
In the present work, we investigate the Lie algebra of the Formanek-Procesi group ${\rm FP}({\rm H})$ with base group $H$ a right-angled Artin group. We show that the Lie algebra $\gr(FP(H))$ has a presentation that is dictated by the group presentation. Moreover we show that ${\rm FP}({\rm H})$ is a Magnus group.
\end{abstract}

%\tableofcontents

%%%%%%%%%%%%%%%%%%%%%%%%%%%%%%%%%%%%
\section{Introduction}
%%%%%%%%%%%%%%%%%%%%%%%%%%%%%%%%%%%%%%

The Formanek-Procesi group was introduced in \cite{fp} as a counterexample to the linearity of Aut$(F_{n})$, the group of automorphisms of a free group $F_{n}$ of rank $n \geq 3$. For any group $G$, the Formanek-Procesi group is the HNN-extension with presentation 
$$
{\rm FP(G)}=\< b, G\times G: b(g,g)b^{-1}=(1,g), \ g\in G\>.
$$ 
It was shown in \cite{fp} that ${\rm FP}(F_{n})$ is not linear.  The group ${\rm FP}(F_{n})$ has been one of the fundamental examples for the non-linearity of groups which is often mentioned as poison group.

A right-angled Artin group (briefly, raag) is a finitely generated group with a standard presentation of the form $H=\< a_1,\ldots, a_k\mid [a_i,a_j]=1~\text{for~ certain~ pairs}~ (i,j)~ \text{with}~ i,j\in\{1,\ldots,k\}\>$. The generating set $\{a_1,\ldots,a_k\}$ is the standard generating set for $H$. The defining graph $\Gamma$ of $H$ is the graph whose vertices correspond to infinite cyclic groups generated by $a_i$ and whose edges connect a pair of vertices $a_i, a_j$ if and only if $[a_i,a_j]=1$. For an extensive survey concerning raags, one may consult \cite{charney} and the references therein. It is obvious that if $H$ is not abelian, then $H$ contains free subgroups and so ${\rm FP(H)}$ is not linear.

In this work, we investigate the associated Lie algebra of ${\rm FP(H)}$ by providing an elimination method for $\gr({\rm FP(H)})$. The elimination method was introduced by Lazard in \cite{laz} for the study of nilpotent groups and $N$-series. Duchamp and Krob \cite{dk} study the free partially commutative Lie algebra and prove a partially commutative version of Lazard's elimination. In the present work, we prove a Lazard's elimination version of the free Lie algebra by decomposing it into a summand of free Lie algebras and ideals that originate by 
the various relators that appear in  the Formanek-Procesi group with base a raag. To achieve this decomposition, besides the classical Lazard's elimination, we also use  results of Labute \cite{labu} and of Duchamp and Krob \cite{dk}. As a result, we get a presentation of the Lie algebra of the Formanek-Procesi group with basis a raag (see Theorem \ref{thm4.12}) and the quotient groups of the consecutive terms of its lower central series are free abelian groups of finite rank (see Corollary \ref{co4.13}).  Furthermore, using a result of Guaschi and Pereiro \cite[Theorem 1.1]{gp}, we show that the Formanek-Procesi group with basis a raag is residually nilpotent (see Corollary \ref{co5.8}). A first attempt in that direction was made in \cite{sev}.

To be more specific, for $n \geq 2$,  let $\Delta =
\{\n\times\n\}\setminus\{(a,a): a\in \n\}$. Let $\theta$ be a non-empty subset of $\Delta$ such that $\theta$ contains both $(a,b)$ and $(b,a)$ for $a,b\in\n$. We call $\theta$ a partial commutation relation on $\n$. Now let $H_1$ and $H_2$ two copies of a raag with presentations
$$
H_i=\<a_{i,1},\ldots, a_{i,n}: [a_{i,k},a_{i,l}]=1,\; \ k > l; (k,l)\in\theta\>, \ i=1,2\>
$$ 
and define ${\rm FP(H)}$ to be the group with generating set $b, a_{1,1},\ldots,a_{1,n},a_{2,1},\ldots, a_{2,n}$ and defining relations
$$
ba_{1,k}a_{2,k}b^{-1}=a_{2,k}, \ \ k\in \n,
$$
$$
[a_{i,r},a_{i,s}]=1,~~ i=1,2;\ \ r > s;  (r,s)\in\theta~~\text{and}
$$
$$
[a_{1,k},a_{2,l}]=1, \ \ k,l\in \n.
$$
Using the first set of the above relation, we may solve for $a_{1,k}$, with $k\in \n$, and replace into the other two sets of relations to get a presentation with generators $b, a_{2,k}$, with $k\in\n$ and relations
$$
[b,a_{2,k},a_{2,l}]=1, [a_{2,r},a_{2,s}]=1,[[b,a_{2,r}],[b,a_{2,s}]]=1~ \mbox{for}\ k,l\in \n; r > s; (r,s)\in\theta.
$$ By renaming $a_{2,i}=f_i$ and $b=t$,  ${\rm FP(H)}$ has a presentation with generating set $\{t, f_{1}, \ldots, f_{n}\}$ and defining relations 
$$ 
[t,f_{k},f_{l}]=[f_{r},f_{s}]=[[t,f_{r}],[t,f_{s}]]=1,  k,l\in\n; r > s; (r,s)\in\theta.
$$

Let $F_{n+1}$ be the free group with free generating set ${\mathcal{F}}_{n+1} = \{f_{1}, \ldots, f_n,f_{n+1}=t\}$, $n \geq 2$, and let $\theta$ be a non-empty partial commutation relation on $\{1, \ldots, n\}$.  Let $M = {\mathcal{R}}^{F_{n+1}}$ be the normal subgroup of ${\mathcal{R}}$ in $F_{n+1}$, where ${\mathcal{R}} = \{[f_{a},f_{b}], [t, f_{i}, f_{j}], [[t,f_{a}], [t,f_{b}]]: i, j \in \{1, \ldots, n\}; a > b; (a,b) \in \theta\}$. 
Notice that ${\rm FP(H)} = F_{n+1}/M$.  
Let ${\rm gr}(F_{n+1})$ be the free Lie algebra (over $\mathbb{Z}$) with a free generating set $\overline{\mathcal{F}}_{n+1} = \{y_{i} = f_{i}F^{\prime}_{n+1}, \overline{t}=tF^{\prime}_{n+1}: i = 1, \ldots, n\}$  associated to $F_{n+1}$.  We order the elements of $\overline{\mathcal{F}}_{n+1}$ as $y_1<\ldots <y_n<\t$. Let also $J$ be the ideal of ${\rm gr}(F_{n+1})$ generated by $\mathcal{R}_{L} = \{[y_{a}, y_{b}], [\t, y_{i}, y_{j}], [[\t, y_{a}], [\t, y_{b}]]: i, j \in \{1, \ldots, n\}; a > b; (a,b) \in \theta\}$.  We use this notation throughout the present work.
 
Our main theorem is the following.

\begin{theorem}\label{mth} ${\rm gr}({\rm FP(H)}) \cong {\rm gr}(F_{n+1})/J$ as Lie algebras in a natural way.  Moreover ${\rm FP(H)}$ is a Magnus group. 
\end{theorem}

The structure of the paper is the following. In section \ref{sect2}, we give some known results concerning Lazard's elimination and decompose ${\rm gr}(F_{n+1})$ for $n\ge 1$ in the spirit of elimination theorem. 
In subsection \ref{subsec2.3}, we analyze the ideal $J$ of ${\rm gr}(F_{n+1})$, $n \geq 1$, as a direct sum of three explicitly described free Lie subalgebras of ${\rm gr}(F_{n+1})$. 
These free Lie subalgebras play an essential role in the proof of Theorem \ref{mth}. 
Let $\theta$ be a non-empty  partial commutation relation on $\{1, \ldots, n\}$, $n \geq 2$, ${\cal{P}} = \{r_{a,b} = [f_{a},f_{b}]: a > b; (a,b) \in \theta\}$ and let $N = {\cal{P}}^{F_{n}}$ be the normal closure of ${\cal{P}}$ in $F_{n}$. 
Let $I$ be the ideal in ${\rm gr}(F_{n})$ generated by the set ${\cal{P}}_{L} = \{[y_{a}, y_{b}]: a > b; (a,b) \in \theta\}$. In section \ref{sect3}, we show that $I$ is a direct summand of ${\rm gr}(F_{n})$ and so ${\rm gr}(F_{n})/I$ is a free $\mathbb{Z}$-module. 
A decomposition of $I$ is given where each summand is a free Lie algebra by explicitly giving a free generating set. 
Furthermore we show that ${\mathcal{L}}(N) = I$. In section \ref{sect4}, we prove ${\rm gr}(F_{n+1}/M) \cong {\rm gr}(F_{n+1})/J$ as Lie algebras in a natural way. 
As a consequence, each ${\rm gr}_{d}({\rm FP(H)})$ is a free abelian group of finite rank. We use results of section \ref{sect3} to show that ${\mathcal{L}}(M) = J$. 
In section \ref{sect5}, we show that ${\rm FP(H)}$ is residually nilpotent  by using a result of Guaschi and Pereiro \cite[Theorem 1.1]{gp}.

%%%%%%%%%%%%%%%%%%%%%%%%%%%%%%%%%%%%%
\section{Preliminaries}\label{sect2}
%%%%%%%%%%%%%%%%%%%%%%%%%%%%%%%%%%%%%

\subsection{The elimination theorem}\label{subsec2.1}

By \emph{a Lie algebra} $L$, we mean a Lie algebra over the ring $\mathbb{Z}$ of integers. Throughout this paper, we use the left-normed convention for Lie commutators. Let $L$ be any Lie algebra. For $x, y \in L$ and $m$ a non-negative integer, we write $[y,~_{m}x] = [y, x, \ldots, x]$ with $m$ factors of $x$. For $m = 0$, we write $[y, ~_{0}x] = y$. Similar notation we use for group commutators. 

Let $V$ be a free $\mathbb{Z}$-module.  We write $L(V)$ for the free Lie algebra which has $V$ as a submodule and every basis of $V$ as a free generating set. 
Thus we may write $L(V) = L({\mathcal{V}})$ for a $\mathbb{Z}$-basis ${\mathcal{V}}$ of $V$. 
By a graded $\mathbb{Z}$-module, we mean a $\mathbb{Z}$-module $V$ with a distinguished $\mathbb{Z}$-module decomposition $V = V_{1} \oplus V_{2} \oplus \cdots$, where each $V_{m}$ is a free $\mathbb{Z}$-module of finite rank. For each $d \geq 1$, let $L^{d}_{\rm grad}(V)$ be the $\mathbb{Z}$-submodule of $L(V)$ spanned by all left-normed Lie commutators $[v_{1}, v_{2}, \ldots, v_{\kappa}]$, with $\kappa \geq 1$, such that, for $i \in \{1,\ldots,\kappa\}$, $v_{i} \in V_{m(i)}$ for some $m(i) \geq 1$ with $m(1) + m(2) + \cdots + m(\kappa) = d$. Then $L(V)$ is a graded $\mathbb{Z}$-module $
L(V) = L^{1}_{\rm grad}(V) \oplus L^{2}_{\rm grad}(V) \oplus \cdots$. 
If $V$ is a free $\mathbb{Z}$-module of finite rank, then we may regard $V$ as a graded $\mathbb{Z}$-module: $V = V \oplus 0 \oplus 0 \oplus \cdots$. Then the $\mathbb{Z}$-module $L^{d}(V)$ is the $d$-th homogeneous component of $L(V)$.  

For non-empty subsets ${\cal B}$ and ${\cal C}$ of a Lie algebra, we write $
{\cal C} \wr {\cal B}^{*} = \{[c, b_{1}, \ldots, b_{k}]: c \in {\cal C}; b_{1}, \ldots, b_{k} \in {\cal B}; k \geq 0\}$.  
The following
result is a version of Lazard's ``Elimination Theorem" (see
\cite[Chapter 2, Section 2.9, Proposition 10]{bour}). 

\begin{theorem}[Elimination]\label{th1}
Let ${\cal A} = {\cal B} \cup {\cal C}$ be the disjoint union of its proper non-empty subsets $\cal B$ and $\cal C$ and consider the free Lie algebra $L({\cal A})$. Then $\cal B$ and ${\cal C} \wr {\cal B}^{*}$ freely generate Lie algebras $L({\cal B})$ and $L({\cal C} \wr {\cal B}^{*})$, and there is a $\mathbb{Z}$-module decomposition $L({\cal A}) = L({\cal B}) \oplus L({\cal C} \wr {\cal B}^{*})$. Furthermore $L({\cal C} \wr {\cal B}^{*})$ is the ideal of $L({\cal A})$ generated by ${\cal C}$. 
\end{theorem}

For non-empty subsets $X, Y$ of $L(\cal{A})$, we write $
S(X, Y) = \{[x,y]: x \in X, y \in Y\}$. We inductively define the subsets $({\cal A}_{m})_{m \geq 1}$ of $L(\cal{A})$ by $
{\cal A}_{1} = {\cal A}$  and, for all $m \geq 2$, $
{\cal A}_{m} = \bigcup_{p, q \geq 1 \atop p+q=m} (S({\cal A}_{p}, {\cal A}_{q}) \setminus \{0\})$.  
An element of $\bigcup_{m \geq 1}{\cal A}_{m}$ is called a simple Lie commutator which is formed by the elements of $\cal{A}$. The following result has been proved in \cite[Lemma II.5]{dk}.

\begin{lemma}\label{le2.2}
Let ${\cal A} = {\cal B} \cup {\cal C}$ be the disjoint union of its proper non-empty subsets $\cal B$ and $\cal C$ and consider the free Lie algebra $L({\cal A})$. Let $u$ be a simple Lie commutator of  $L({\cal A})$. Then either $u \in L({\cal B})$ or $u \in L({\cal C} \wr {\cal B}^{*})$. 
\end{lemma}

%%%%%%%%%%%%%%%%%%%%%%%%%%%%%%%%%%%%%
\subsection{Decomposition of ${\rm gr}(F_{n+1})$}\label{sec2.2}
%%%%%%%%%%%%%%%%%%%%%%%%%%%%%%%%%%%%%

Let $G$ be any group. For subgroups $X$ and $Y$ of $G$, we write $[X,Y]$ for the subgroup of $G$ generated by the commutators $[x,y] = x^{-1}y^{-1}xy$ with $x \in X$ and $y \in Y$. 
For $c \geq 1$, let $\gamma_{c}(G)$ be the $c$-th term of the lower central series of $G$. 
That is, $\gamma_{1}(G) = G$ and, for $c \geq 2$, $\gamma_{c}(G) = [\gamma_{c-1}(G),G]$. 
We write $G^{\prime} = \gamma_{2}(G)$ for the commutator subgroup of $G$. 
For subgroups $X_{1}, \ldots, X_{m}$ of $G$, with $m \geq 2$, we write $[X_{1}, X_{2}, \ldots, X_{m}] = [[X_{1}, \ldots, X_{m-1}],X_{m}]$ for the subgroup of $G$ generated by all group commutators $[u,x]$ where $u \in [X_{1}, \ldots, X_{m-1}]$ and $x \in X_{m}$. 
Let $Y$ be a subgroup of $G$. 
For $\kappa \geq 1$, we write $[Y,~_{\kappa}G] = [Y, G, \ldots, G]$ with $\kappa$ factors of $G$.  
Then $[Y,\prescript{}{\k}{G}]$ is a normal subgroup of $G$ with a generating set $\{[y, a_{1}, \ldots, a_{\k}]: y \in Y; a_{1}, \ldots, a_{\k} \in G\}$. 
For a non-empty subset $\mathcal{X}$ of $G \setminus \{1\}$ we write  $\mathcal{X}^{G}$ for the normal closure of $\mathcal{X}$ in $G$. 
In particular, if $X$ is the subgroup of $G$ generated by $\mathcal{X}$, then  $\mathcal{X}^{G} = X~[X,G]$. 
The associated graded abelian group ${\rm gr}(G) = \bigoplus_{c \geq 1}{\rm gr}_{c}(G)$, where ${\rm gr}_{c}(G) = \gamma_{c}(G)/\gamma_{c+1}(G)$, has the structure of a graded Lie algebra over $\mathbb{Z}$, the Lie bracket operation in ${\rm gr}(G)$ being induced by the commutator operation in $G$ (see, for example, \cite{bour}, \cite{laz}, \cite{mks}). 

The proof  of the following result is elementary.

\begin{lemma}\label{le3.3}
Let $G$ be a group and $Y$ be a subgroup of $G$. Then,  
for $m, r \ge 1$, $[Y, \gamma_{m}(G)] \leq [Y, ~_{m}G]$ and $[\gamma_{r}(G), Y, \gamma_{m}(G)] \leq [[\gamma_{r}(G), Y],~_{m}G] \leq [Y, ~_{m+r}G]$. 
\end{lemma}

The Lie algebra ${\rm gr}(F_{n+1})$ is a free Lie algebra with a free generating set $\overline{\mathcal{F}}_{n+1} = \{\t, y_{i} : i = 1, \ldots, n\}$.  For $\k\in\{2,\ldots, n+1\}$ let $F_{\k}$ be the free group generated by ${\mathcal F}_{\k}$.  
For $c \geq 1$, ${\cal L}_{{\rm gr},c}(F_{\kappa}) = \gamma_{c}(F_{\kappa})\gamma_{c+1}(F_{n+1})/\gamma_{c+1}(F_{n+1})$. 
Form the restricted direct sum $
{\cal L}_{{\rm gr}}(F_{\kappa}) = \bigoplus_{c \geq 1}{\cal L}_{{\rm gr},c}(F_{\kappa})$                
of the abelian groups ${\cal L}_{{\rm gr},c}(F_{\kappa})$.  Since $[\gamma_{\mu}(F_{\kappa}), \gamma_{\nu}(F_{\kappa})] \subseteq \gamma_{\mu + \nu}(F_{\kappa})$ for all $\mu, \nu \ge 1$, ${\cal L}_{{\rm gr}}(F_{\kappa})$ is a Lie subalgebra of ${\rm gr}(F_{n+1})$. 
Notice that ${\mathcal{L}}_{{\rm gr}}(F_{n+1}) = {\rm gr}(F_{n+1})$. Since $\gamma_{c+1}(F_{\kappa}) = \gamma_{c}(F_{\kappa}) \cap \gamma_{c+1}(F_{n+1})$ for all $c \geq 1$, we have ${\mathcal{L}}_{{\rm gr}}(F_{\kappa}) \cong {\rm gr}(F_{\kappa})$ as Lie algebras in a natural way. 
Hence ${\mathcal{L}}_{{\rm gr}}(F_{\kappa})$ is a free Lie algebra on $\overline{\mathcal{F}}_{\kappa}$. 
By Theorem \ref{th1}, $
{\mathcal{L}}_{{\rm gr}}(F_{\kappa}) = L(\overline{\mathcal{F}}_{
\kappa-1}) \oplus L(\{y_{\kappa}\} \wr (\overline{\mathcal{F}}_{\kappa-1})^{*})$. The free Lie algebra $L(\{y_{\kappa}\} \wr (\overline{\mathcal{F}}_{\kappa-1})^{*})$ is the ideal of ${\mathcal{L}}_{{\rm gr}}(F_{\kappa})$ generated by the element $y_{\kappa}$. 

\vskip .120 in 

The proof of the following result is straightforward.

\begin{lemma}\label{le2.5a}
For $n \geq 1$, $F_{n+1}$ is the semidirect product of $N_{t}$ by $F_{n}$ where $N_{t}$ is the normal closure of $\{t\}$ in $F_{n+1}$. 
\end{lemma}

\vskip .120 in

Fix $\kappa \in \{2, \ldots, n+1\}$. For $c \geq 1$, $N_{f_{\kappa},c} = N_{f_{\kappa}} \cap \gamma_{c}(F_{\kappa})$, where $N_{f_{\kappa}}$ is the normal closure of $\{f_{\kappa}\}$ in $F_{\kappa}$, and ${\rm I}_{c}(N_{f_{\kappa}}) = N_{f_{\kappa},c} \gamma_{c+1}(F_{n+1})/\gamma_{c+1}(F_{n+1})$. 
Notice that ${\rm I}_{c}(N_{f_{\kappa}}) \leq {\mathcal{L}}_{{\rm gr},c}(F_{\kappa})$ for all $c \geq 1$. Let ${\cal{L}}(N_{f_{\kappa}}) = \bigoplus_{c \geq 1}{\rm I}_{c}(N_{f_{\kappa}})$. 
Since $N_{f_{\kappa}}$ is normal in $F_{\kappa}$, we have ${\cal{L}}(N_{f_{\kappa}})$ is an ideal of ${\mathcal{L}}_{\rm gr}(F_{\kappa})$.

\begin{proposition}\label{pro2.6} For $n+1 \geq 2$ and $\kappa \in \{2, \ldots, n+1\}$, ${\cal L}_{{\rm gr}}(F_{\kappa}) = {\cal L}_{{\rm gr}}(F_{\kappa-1}) \oplus {\cal L}(N_{f_{\kappa}})$. In particular ${\cal L}(N_{f_{\kappa}}) = L(\{y_{\kappa}\} \wr (\overline{\mathcal{F}}_{\kappa-1})^{*})$.  

\end{proposition}

\pf By Lemma \ref{le2.5a} we obtain the following (split) short exact sequence
$$
\{1\} \rightarrow N_{f_{\kappa}} \overset{\varepsilon}{\rightarrow} F_{\kappa} \overset{p}{\rightarrow} F_{\kappa-1} \rightarrow \{1\}, \eqno(2.1)
$$
where $p$ is the epimorphism from $F_{\kappa}$ onto $F_{\kappa-1}$ satisfying the conditions $p(f_{i}) = f_{i}$, $i \in \{1, \ldots, \kappa-1\}$, and $p(f_{\kappa}) = 1$ and $\varepsilon$ is the inclusion map from $N_{f_{\kappa}}$ into $F_{\kappa}$. 
Let $Y_{\kappa}$ be the subgroup of $F_{\kappa}$ generated by $f_{\kappa}$. 
Since $N_{f_{\kappa}} = Y_{\kappa} [Y_{\kappa},F_{\kappa}]$ and $[Y_{\kappa},F_{\kappa}] \subseteq F^{\prime}_{\kappa}$, we have ${\rm I}_{1}(N_{f_{\kappa}}) = \langle y_{\kappa} \rangle$. 
Clearly $[Y_{\kappa}, ~_{(c-1)}F_{\kappa}] \subseteq N_{f_{\kappa},c} \subseteq N_{f_{\kappa}}$ for all $c \geq 2$. 
Since $N_{f_{\kappa}} \cap F_{\kappa-1} = \{1\}$, we get $[Y_{\kappa}, ~_{(c-1)}F_{\kappa}] \cap \gamma_{c}(F_{\kappa-1}) = \{1\}$ for all $c \geq 2$. 
For $r \geq 1$, let $G_{1} = F_{\kappa}$ and, for $r \geq 2$, $G_{r} = [Y_{\kappa}, ~_{(r-1)}F_{\kappa}] \rtimes \gamma_{r}(F_{\kappa-1})$. Notice that $G_{1} \supseteq G_{2} \supseteq \cdots \supseteq G_{r} \supseteq \cdots$. Using the group commutator identity $[ab,cd] = [a,d]^{b} [b,d] [a,c]^{bd}[b,c]^{d}$ and Lemma \ref{le3.3}, we obtain $(G_{r})_{r \geq 1}$ is a strongly central series of $F_{\kappa}$ and so $G_{r} \supseteq \gamma_{r}(F_{\kappa})$ for all $r$. 
Obviously $G_{r} = \gamma_{r}(F_{\kappa})$ for all $r$. Since $[Y_{\kappa}, ~_{(c-1)}F_{\kappa}] \subseteq N_{f_{\kappa}}$ for all $c \geq 2$, we have, by the modular law, $N_{f_{\kappa}} \cap G_{c} = [Y_{\kappa}, ~_{(c-1)}F_{\kappa}](N_{f_{\kappa}} \cap \gamma_{c}(F_{\kappa-1}))$. Since $N_{f_{\kappa}} \cap \gamma_{c}(F_{\kappa-1}) = \{1\}$ and $G_{c} = \gamma_{c}(F_{\kappa})$ for all $c$, we get $[Y_{\kappa}, ~_{(c-1)}F_{\kappa}] = N_{f_{\kappa},c}$ for all $c \geq 2$. Hence ${\rm I}_{d}(N_{f_{\kappa}}) = [Y_{\kappa}, ~_{(c-1)}F_{\kappa}]\gamma_{c+1}(F_{n+1})/\gamma_{c+1}(F_{n+1})$ for all $c \geq 2$.    

As before,  ${\cal{L}}_{{\rm gr}}(F_{\kappa-1}) \cong {\rm gr}(F_{\kappa-1})$ as Lie algebras in a natural way. 
By definition of ${\cal{L}}_{{\rm gr}}(F_{\kappa-1})$ and since $L(\overline{\mathcal{F}}_{\kappa-1})$ is a free Lie algebra on $\overline{\mathcal{F}}_{\kappa-1}$, we obtain ${\cal{L}}_{{\rm gr}}(F_{\kappa-1}) = L(\overline{\mathcal{F}}_{\kappa-1})$. 
Let $p_{L}$ be the Lie algebra epimorphism from ${\mathcal{L}}_{{\rm gr}}(F_{\kappa})$ onto ${\cal{L}}_{{\rm gr}}(F_{\kappa-1})$ satisfying the conditions $p_{L}(y_{i}) = y_{i}$ (for $i \in \{1, \ldots, \kappa-1\}$) and $p_{L}(y_{\kappa}) = 0$. 
Thus ${\rm Ker}p_{L} = L(\{y_{\kappa}\} \wr (\overline{\mathcal{F}}_{\kappa-1})^{*})$. Hence we obtain the following short exact sequence of Lie algebras
$$
\{0\} \rightarrow L(\{y_{\kappa}\} \wr (\overline{\mathcal{F}}_{\kappa-1})^{*}) \overset{\varepsilon_{L}}{\longrightarrow} {\cal{L}}_{{\rm gr}}(F_{\kappa}) \overset{p_{L}}{\longrightarrow} {\cal{L}}_{{\rm gr}}(F_{\kappa-1}) \rightarrow \{0\}, \eqno(2.2)
$$
where $\varepsilon_{L}$ is the inclusion map from $L(\{y_{\kappa}\} \wr (\overline{\mathcal{F}}_{\kappa-1})^{*})$ into ${\cal{L}}_{{\rm gr}}(F_{\kappa})$. Clearly, for all $c \geq 1$, we get by Eq.$(2.2)$ the following short exact sequence of $\mathbb{Z}$-modules 
$$
\{0\} \rightarrow L^{c}_{\rm grad}(\{y_{\kappa}\} \wr (\overline{\mathcal{F}}_{\kappa-1})^{*}) \overset{\varepsilon_{L,c}}{\longrightarrow} {\mathcal{L}}_{{\rm gr},c}(F_{\kappa}) \overset{p_{L,c}}{\longrightarrow} {\mathcal{L}}_{{\rm gr},c}(F_{\kappa-1}) \rightarrow \{0\}, \eqno(2.3)
$$
where $p_{L,c}$ and $\varepsilon_{L,c}$ are the induced group homomorphisms by $p_{L}$ and $\varepsilon_{L}$, respectively. 

For $c \geq 1$, let $p_{c}$ be the group epimorphism from ${\mathcal{L}}_{{\rm gr},c}(F_{\kappa})$ onto ${\mathcal{L}}_{{\rm gr},c}(F_{\kappa-1})$ induced by $p$ (Eq.  $(2.1)$). Further we write $\varepsilon_{c}$ for the group monomorphism from ${\rm I}_{c}(N_{f_{\kappa}})$ into ${\mathcal{L}}_{{\rm gr},c}(F_{\kappa})$ induced by $\varepsilon$ (Eq.  $(2.1)$). 
By definition of $(N_{f_{\kappa},c})_{c \geq 1}$, we have ${\rm Im}\varepsilon_{c} = {\rm Ker}p_{c} = {\rm I}_{c}(N_{f_{\kappa}})$ for all $c$, that is, we obtain the following short exact sequence of $\mathbb{Z}$-modules for all $c$
$$
\{0\} \rightarrow {\rm I}_{c}(N_{f_{\kappa}}) \overset{\varepsilon_{c}}{\longrightarrow} {\mathcal{L}}_{{\rm gr},c}(F_{\kappa}) \overset{p_{c}}{\longrightarrow} {\mathcal{L}}_{{\rm gr},c}(F_{\kappa-1}) \rightarrow \{0\}. \eqno(2.4)
$$
Since $L^{c}_{\rm grad}(\{y_{\kappa}\} \wr (\overline{\mathcal{F}}_{\kappa-1})^{*}) \subseteq {\rm I}_{c}(N_{f_{\kappa}})$, we have by Eqs $(2.4)$ and $(2.3)$ that ${\rm I}_{c}(N_{f_{\kappa}}) = L^{c}_{\rm grad}(\{y_{\kappa}\} \wr (\overline{\mathcal{F}}_{\kappa-1})^{*})$ for all $c$. 
Therefore ${\cal{L}}(N_{f_{\kappa}}) = L(\{y_{\kappa}\} \wr (\overline{\mathcal{F}}_{\kappa-1})^{*})$. 
By Theorem \ref{th1} and  since ${\cal{L}}_{{\rm gr}}(F_{\kappa-1}) = L(\overline{\mathcal{F}}_{\kappa-1})$, we get ${\cal L}_{{\rm gr}}(F_{\kappa}) = {\cal L}_{{\rm gr}}(F_{\kappa-1}) \oplus {\cal L}(N_{f_{\kappa}})$ for all $\kappa \in \{2, \ldots, n+1\}$.~\qed

%%%%%%%%%%%%%%%%%%%%%%%%%%%%%%%%%%%%
\subsection{Analysis of $J$}\label{subsec2.3}
%%%%%%%%%%%%%%%%%%%%%%%%%%%%%%%%%%%%

%For $n \geq 2$, let $\theta$ be a non-empty partial commutation relation on the set $\{1, \ldots, n+1\}$ with $(n+1,i), (i, n+1) \notin \theta$ for all $i \in \{1, \ldots, n\}$. Let $J$ be the ideal of ${\rm gr}(F_{n+1})$ generated by $\mathcal{R}_{L} = \{[\overline{f}_{a}, \overline{f}_{b}], [\overline{f}_{n+1}, \overline{f}_{i}, \overline{f}_{j}], [[\overline{f}_{n+1}, \overline{f}_{a}], [\overline{f}_{n+1}, \overline{f}_{b}]]: i, j \in \{1, \ldots, n\}; a > b; (a,b)$ $\in \theta\}$. Write $t = \overline{f}_{n+1}$ and $y_{k} = \overline{f}_{k}$ for $k \in \{1, \ldots, n\}$. Thus $\overline{\mathcal{F}}_{n+1} = \{y_{1}, \ldots, y_{n}, t\}$, and order its elements as $y_{1} < y_{2} <  \cdots < y_{n} < t$. Notice that $\overline{\mathcal{F}}_{n} = \{y_{1}, \ldots, y_{n}\}$. 

By Proposition \ref{pro2.6}  
$$
{\rm gr}(F_{n+1}) = \mathcal{L}_{{\rm gr}}(F_{n}) \oplus L(\{\t\} \wr (\overline{\mathcal{F}}_{n})^{*}). \eqno(2.5)
$$
Let $\Omega = \{\t, [\t,y_{i}]: i \in \n\}$ and $\mathcal{E} = (\{\t\} \wr (\overline{\mathcal{F}}_{n})^{*}) \setminus \Omega$. By Theorem \ref{th1}, we have 
$$
L(\{\t\} \wr (\overline{\mathcal{F}}_{n})^{*}) = L(\Omega \cup {\cal{E}}) = L(\Omega) \oplus L({\cal E} \wr \Omega^{*}). \eqno(2.6)
$$
The free Lie algebra $L({\cal E} \wr \Omega^{*})$ is the ideal in $L(\{\t\} \wr (\overline{\mathcal{F}}_{n})^{*})$ generated by $\cal E$. Since ${\cal E} \subseteq J$ and $J$ is an ideal in ${\rm  gr}(F_{n+1})$, we have $L({\cal E} \wr \Omega^{*}) \subseteq J$. By Eq. $(2.6)$ and the modular law, we get 
$$
J \cap L(\{\t\} \wr (\overline{\mathcal{F}}_{n})^{*}) = (J \cap L(\Omega)) \oplus L({\cal E} \wr \Omega^{*}). \eqno(2.7)
$$
By the decomposition of $L(\{\t\} \wr (\overline{\mathcal{F}}_{n})^{*})$ (Eq. $(2.6)$), the equation $(2.5)$ becomes
$$
{\rm gr}(F_{n+1}) = \mathcal{L}_{{\rm gr}}(F_{n}) \oplus L(\Omega) \oplus L({\cal E} \wr \Omega^{*}). \eqno(2.8)
$$

\begin{lemma}\label{le3.5}
The ideal $J \cap L(\{\t\} \wr (\overline{\mathcal{F}}_{n})^{*})$ of ${\rm gr}(F_{n+1})$ is generated by $\{[[\t, y_{a}], [\t, y_{b}]],$ $[\t, y_{i}, y_{j}]: i, j \in \n; a > b; (a,b) \in \theta\}$.
\end{lemma}

\pf Let $J_{1}$ and $J_{2}$ be the ideals of ${\rm gr}(F_{n+1})$ generated by $
{\cal{J}}_{1} = \{[y_{a},y_{b}]: a > b; (a,b) \in \theta\}$ 
and $
{\cal{J}}_{2} = \{[\t, y_{i}, y_{j}], [[\t,y_{a}],[\t,y_{b}]]: i,j \in \n; a > b; (a,b) \in \theta\}$, 
respectively. 
Clearly $J = J_{1} + J_{2}$. Since $L(\{\t\} \wr (\overline{\mathcal{F}}_{n})^{*})$ is an ideal in ${\rm gr}(F_{n+1})$ and ${\cal{J}}_{2} \subseteq L(\{\t\} \wr (\overline{\mathcal{F}}_{n})^{*})$, we conclude that $J_{2} \subseteq L(\{\t\} \wr (\overline{\mathcal{F}}_{n})^{*})$ and so $J_{2} \subseteq J \cap L(\{\t\} \wr (\overline{\mathcal{F}}_{n})^{*})$. 
But
$$
\begin{array}{rll}
J \cap L(\{\t\} \wr (\overline{\mathcal{F}}_{n})^{*}) & = & (J_{1} + J_{2}) \cap L(\{\t\} \wr (\overline{\mathcal{F}}_{n})^{*}) \\
({\rm modular~~law})& = & (J_{1} \cap L(\{\t\} \wr (\overline{\mathcal{F}}_{n})^{*})) + J_{2}.
\end{array}
$$ 
In order to show that $J \cap L(\{\t\} \wr (\overline{\mathcal{F}}_{n})^{*}) \subseteq J_{2}$, it is enough to prove that $J_{1} \cap L(\{\t\} \wr (\overline{\mathcal{F}}_{n})^{*}) \subseteq J_{2}$. By Lemma \ref{le2.2} (and Eq. $(2.5)$), we have  $
J_{1} = (J_{1} \cap \mathcal{L}_{{\rm gr}}(F_{n})) \oplus (J_{1} \cap L(\{\t\} \wr (\overline{\mathcal{F}}_{n})^{*})$. Clearly $J_{1} \cap L(\{\t\} \wr (\overline{\mathcal{F}}_{n})^{*})$ is an ideal of ${\rm gr}(F_{n+1})$. 
By a result of Duchamp and Krob \cite[Theorem II.7]{dk}, we have $J_{1} \cap \mathcal{L}_{{\rm gr}}(F_{n})$ is the ideal of $\mathcal{L}_{{\rm gr}}(F_{n})$ generated by ${\cal{J}}_{1}$ and $J_{1} \cap L(\{\t\} \wr (\overline{\mathcal{F}}_{n})^{*})$ is the ideal of $L(\{\t\} \wr (\overline{\mathcal{F}}_{n})^{*})$ generated by ${\cal{T}} = \{[\t,z_{1}, \ldots,z_{\mu},[y_{a},y_{b}],z^{\prime}_{1}, \ldots,z^{\prime}_{\nu}]: z_{1} \cdots z_{\mu}, z^{\prime}_{1} \cdots z^{\prime}_{\nu} \in (\overline{\mathcal{F}}_{n})^{*}, (a,b) \in \theta\}$.  
(In the statement of Theorem II.7 of \cite{dk}, the set $T_{1} = \emptyset$, since in our case $(n+1,i), (i,n+1) \notin \theta$ for all $i \in \{1, \ldots, n\}$.) 
Using the Jacobi identity in the form $[x,[y,z]] = [x,y,z] - [x,z,y]$, we get ${\cal{T}} \subseteq J_{2}$. 
Since $J_{2}$ is an ideal in ${\rm gr}(F_{n+1})$, we get $J_{1} \cap L(\{\t\} \wr (\overline{\mathcal{F}}_{n})^{*}) \subseteq J_{2}$ and so $J \cap L(\{\t\} \wr (\overline{\mathcal{F}}_{n})^{*}) = J_{2}$. 
That is, $J \cap L(\{\t\} \wr (\overline{\mathcal{F}}_{n})^{*})$ is the ideal in ${\rm gr}(F_{n+1})$ generated by ${\cal J}_{2}$.~\qed

As before,  by $F_n$ we denote the subgroup of $F_{n+1}$ generated by ${\mathcal F}_n = \{f_1,\ldots,f_n\}$.

\begin{proposition}\label{pro3.6} 
\begin{enumerate} 

\item $J = (J \cap \mathcal{L}_{{\rm gr}}(F_{n})) \oplus (J \cap L(\Omega)) \oplus L({\mathcal{E}} \wr \Omega^{*})$.

\item The ideal $J \cap \mathcal{L}_{{\rm gr}}(F_{n})$ of $\mathcal{L}_{{\rm gr}}(F_{n})$ is generated by ${\mathcal{R}}_{2,L} = \{[y_{a},y_{b}]: a > b; (a, b) \in \theta\}$.

\item $J \cap L(\Omega)$ is the ideal of $L(\Omega)$ generated by ${\mathcal{R}}_{4,L} = \{[[\t, y_{a}], [\t, y_{b}]]: a > b; (a,b) \in \theta\}$. 

\item $L({\cal E} \wr \Omega^{*})$ is the ideal of ${\rm gr}(F_{n+1})$ generated by ${\mathcal{R}}_{3,L} = \{[\t,y_{i},y_{j}]: i, j \in \{1, \ldots, n\}\}$.
\end{enumerate} 
\end{proposition}

\pf \begin{enumerate}

\item By Eq. $(2.5)$, the definition of $J$ and  Lemma \ref{le2.2}, we have $$
J = (J \cap \mathcal{L}_{{\rm gr}}(F_{n})) \oplus (J \cap L(\{\t\} \wr (\overline{\mathcal{F}}_{n})^{*})). \eqno(2.9)
$$   
Using Eq. $(2.7)$, we obtain the required decomposition of $J$. 

\item Since $J_{1} \subseteq J$, we have $J_{1} \cap {\mathcal{L}}_{\rm gr}(F_{n}) \subseteq J \cap {\mathcal{L}}_{\rm gr}(F_{n})$. 
By Eq. $(2.9)$, Lemma \ref{le3.5} and the definition of $J$, we get $J \cap {\mathcal{L}}_{\rm gr}(F_{n}) \subseteq J_{1} \cap {\mathcal{L}}_{\rm gr}(F_{n})$ and so $J_{1} \cap {\mathcal{L}}_{\rm gr}(F_{n}) = J \cap {\mathcal{L}}_{\rm gr}(F_{n})$. 
By the proof of Lemma \ref{le3.5}, we obtain the desired result.   

\item By Lemma \ref{le3.5}, we have $J \cap L(\{\t\} \wr (\overline{\mathcal{F}}_{n})^{*}) \subseteq \bigoplus_{m \geq 3}{\rm gr}_{m}(F_{n+1})$.  
Let $\phi$ be the natural map from $L(\{\t\} \wr (\overline{\mathcal{F}}_{n})^{*})$ onto $L(\Omega)$. 
Thus $\phi(\t) = \t$, $\phi([\t,y_{i}]) = [\t,y_{i}]$, $i \in \n$, and $\phi([\t,b_{1}, \ldots, b_{k}]) = 0$ for all $b_{i} \in \overline{{\mathcal F}}_n=\{y_1,\ldots, y_n\}$, with $k \geq 2$. 
Since $\phi$ is a Lie algebra epimorphism, we have $\phi(J \cap L(\{\t\} \wr (\overline{\mathcal{F}}_{n})^{*}))$ is the ideal of $L(\Omega)$ generated by the set $\{[[\t,y_{a}],[\t,y_{b}]]: a>b; (a,b) \in \theta\}$. 
It is clear enough that $\phi(J \cap L(\{\t\} \wr (\overline{\mathcal{F}}_{n})^{*})) \subseteq J \cap L(\Omega)$. Let $y \in J \cap L(\Omega)$. 
There exists $x \in L(\{\t\} \wr (\overline{\mathcal{F}}_{n})^{*})$ such that $\phi(x) = y$. 
By Eq. $(2.6)$, we get $x = u + v$, where $u \in L(\Omega)$ and $v \in L({\cal E} \wr \Omega^{*})$. 
Since $\phi({\cal E} \wr \Omega^{*}) = 0$ and $L({\cal E} \wr \Omega^{*})$ is a free Lie algebra on ${\cal E} \wr \Omega^{*}$, we have $\phi(x) = \phi(u) + \phi(v) = u = y$ and so $u \in J \cap L(\Omega)$. 
By Eq. $(2.7)$, we conclude that $x \in J \cap L(\{\t\} \wr (\overline{\mathcal{F}}_{n})^{*})$. Hence $y \in \phi(L(\{\t\} \wr (\overline{\mathcal{F}}_{n})^{*}))$ and so $
\phi(J \cap L(\{\t\} \wr (\overline{\mathcal{F}}_{n})^{*})) = J \cap L(\Omega)$. 
Therefore the ideal $J \cap L(\Omega)$ of $L(\Omega)$ is generated by the set $\{[[\t,y_{a}],[\t,y_{b}]]: a > b; (a,b) \in \theta\}$.

\item Throughout the proof, we write $S = L({\cal E} \wr \Omega^{*})$. To show that $S$ is an ideal of ${\rm gr}(F_{n+1})$, it is enough to show that $[u,v] \in S$ for all $u \in S$ and $v \in {\rm gr}(F_{n+1})$. By Eq. $(2.5)$ and since $S$ is an ideal of $L(\{\t\} \wr (\overline{\mathcal{F}}_{n})^{*})$, it is enough to show that $[u,v] \in S$ for all $u \in S$ and $v \in \mathcal{L}_{\gr}(F_{n})$. Since $\mathcal{L}_{{\rm gr}}(F_{n})$ is generated by the set $\overline{\mathcal{F}}_{n}$, every element $v$ of $\mathcal{L}_{{\rm gr}}(F_{n})$ is a $\mathbb{Z}$-linear combination of Lie commutators of the form $[y_{i_{1}}, \ldots, y_{i_{\kappa}}]$, $\kappa \geq 1$, $y_{i_{1}}, \ldots, y_{i_{\kappa}} \in \overline{\mathcal{F}}_{n}$. Using the Jacobi identity in the form $
[x, [y,z]] = [x,y,z] - [x,z,y]$, we write each Lie commutator of the form $[u,[y_{i_{1},1}, \ldots, y_{i_{k_{1}},1}], \ldots, [y_{i_{1},\mu}, \ldots, y_{i_{k_{\mu}},\mu}]]$ as a $\mathbb{Z}$-linear combination of Lie commutators of the form $[u,y_{j_{1},1}, \ldots, y_{j_{k_{1}},1}, \ldots, y_{j_{1},\mu}, \ldots, y_{j_{k_{\mu}},\mu}]$. By the multi-linearity of Lie commutators, it is enough to show that $[u, y_{i_{1}}, \ldots, y_{i_{\kappa}}] \in S$ for all $u \in S$ and $y_{i_{1}}, \ldots, y_{i_{\kappa}} \in \overline{{\mathcal F}}_n$, $\kappa \geq 1$. For a positive integer $\alpha$, let $[{\cal{E}},~_{\alpha}\Omega] = [{\cal E}, \Omega, \ldots, \Omega] = \{[v, z_{1}, \ldots, z_{\alpha}]: v \in {\cal{E}}, z_{1}, \ldots, z_{\alpha} \in \Omega\}$. 
Since $
{\cal E} \wr \Omega^{*} = {\cal E} \cup (\bigcup_{\alpha \geq 1}[{\cal{E}},~_{\alpha}\Omega])$ 
and $S$ is a free Lie algebra on ${\cal E} \wr \Omega^{*}$, it is enough to show that $[w, y_{i_{1}}, \ldots, y_{i_{\kappa}}] \in S$ for all $w \in {\cal E} \wr \Omega^{*}$ and $y_{i_{1}}, \ldots, y_{i_{\kappa}} \in \overline{\mathcal{F}}_{n}$, $\kappa \geq 1$. 
Since ${\cal E} \wr \Omega^{*}$ is a free generating set of $S$, we may assume that $w \in [{\cal{E}},~_{\alpha}\Omega]$ with $\alpha \geq 1$. 
Since $S$ is an ideal of $L(\{\t\} \wr (\overline{\mathcal{F}}_{n})^{*})$ and using the Jacobi identity in the form $
[x,y,z] = [x,z,y] + [x,[y,z]]$  on the elements $[v, z_{1}, \ldots, z_{\alpha}]$ with $v \in {\cal{E}}$ and $z_{1}, \ldots, z_{\alpha} \in \Omega$, we may take $w$ to have the form $
w = [u, ~_{m}\t, ~_{\mu_{1}}[\t,y_{1}], \ldots, ~_{\mu_{n}}[\t,y_{n}]]$, 
where $u \in {\cal E}$, $m, \mu_{1}, \ldots, \mu_{n} \geq 0$ and $m + \mu_{1} + \cdots + \mu_{n} \geq 1$. 
Using the aforementioned Jacobi identity and since $S$ is an ideal of $L(\{\t\} \wr (\overline{\mathcal{F}}_{n})^{*})$, we get $
[w, y_{i_{1}}, \ldots, y_{i_{\kappa}}] = [u, y_{i_{1}}, \ldots, y_{i_{\kappa}}, \prescript{}{m}{\t}, \prescript{}{\mu_{1}}{[\t,y_{1}]}, \ldots, \prescript{}{\mu_{n}}[{\t,y_{n}]]} + v^{\prime}$  
where $v^{\prime} \in S$. Since $
[u, y_{i_{1}}, \ldots, y_{i_{\kappa}}, \prescript{}{m}{\t}, \prescript{}{\mu_{1}}{[\t,y_{1}]}, \ldots, \prescript{}{\mu_{n}}{[\t,y_{n}]]} \in S$,  
we conclude that $[w, y_{i_{1}}, \ldots, y_{i_{\kappa}}] \in S$. Therefore $L({\cal E} \wr \Omega^{*})$ is an ideal of ${\rm gr}(F_{n+1})$. \qed 

\end{enumerate}

\section{Free partially commutative Lie algebra}\label{sect3}

As before let $F_n=\<f_1,\ldots, f_n\>$ be a free group and  $\theta$ be a non-empty  partial commutation relation on $\{1, \ldots, n\}$, with $n \geq 2$. 
Let also ${\cal{P}} = \{r_{a,b} = [f_{a},f_{b}]: a > b; (a,b) \in \theta\}$ and let $N = {\cal{P}}^{F_{n}}$ be the normal closure of ${\cal{P}}$ in $F_{n}$. Thus $N \subseteq F^{\prime}_{n}$. For $d \geq 1$, define ${\rm I}_{d}(N) = (N \cap \gamma_{d}(F_{n}))\gamma_{d+1}(F_{n})/\gamma_{d+1}(F_{n}) \leq {\rm gr}_{d}(F_{n})$ and note that ${\rm I}_{1}(N) = \{0\}$. 
Form the (restricted) direct sum ${\cal L}(N) = \bigoplus_{d \geq 2}{\rm I}_{d}(N)$ of the abelian groups ${\rm I}_{d}(N)$. 
Since $N$ is a normal subgroup of $F_{n}$, the Lie subalgebra ${\cal L}(N)$ is an ideal of ${\rm gr}(F_{n})$. Let $I$ be the ideal in ${\rm gr}(F_{n})$ generated by ${\cal{P}}_{L} = \{[y_{a}, y_{b}]: a > b; (a,b) \in \theta\}$. 
Since ${\cal{P}}_{L} \subseteq {\rm I}_{2}(N)$ and ${\cal{L}}(N)$ is an ideal of ${\rm gr}(F_{n})$, we get $I \subseteq {\cal{L}}(N)$. 

Our aims in this section are to show $I$ is a direct summand of ${\rm gr}(F_{n})$ and ${\mathcal{L}}(N) = I$. Proving that $I$ is a direct summand of ${\rm gr}(F_{n})$, we have ${\rm gr}(F_{n})/I$ is a free $\mathbb{Z}$-module. Since $I$ is a homogeneous ideal, we get, by a result of Witt (see, for example, \cite[Section 2.4, Theorem 5. See, also, Theorem 3]{baht}), $I$ is a free Lie algebra. We present a decomposition of $I$ where each of its summands is a free Lie algebra by explicitly giving a free generating set.

If $\theta = \emptyset$, then $I = \{0\}$ and so there is nothing to prove.  

So, without loss of generality, we may assume that $(2,1)\in\theta$. It is convenient throughout this section to write $\overline{\mathcal{F}}_{n} = \mathcal{Y}$. 
%and $\overline{f}_{i} = y_{i}$ for $i \in \{1, \ldots, m\}$. Order the elements of ${\cal Y}$ as $y_1<y_2<\cdots<y_m$. 
For $i \in \{1, \ldots, n\}$, let ${\cal Y}_{i} = \{y_{1}, y_{2}, \ldots, y_{i}\}$. Thus ${\cal Y}_{n} = {\cal Y}$. Furthermore, for $i \in \{2, \ldots, n\}$, we write ${\cal I}_{{\cal Y}_{i}} = \{[y_a,y_b], [y_{a},y_{b}]: y_a, y_b \in {\cal Y}_{i}; (a,b) \in \theta\}$ and  $I_{{\cal Y}_{i}}$ for the ideal of ${\mathcal{L}}_{\rm gr}(F_{i})$ generated by the set ${\cal I}_{{\cal Y}_{i}}$. In particular $I = I_{{\cal Y}_{n}}$. Since ${\cal Y}_{2} \subset {\cal Y}_{3} \subset \cdots \subset {\cal Y}_{n} = {\cal Y}$, we have $
{\cal I}_{{\cal Y}_{2}} \subseteq {\cal I}_{{\cal Y}_{3}} \subseteq \cdots \subseteq {\cal I}_{{\cal Y}_{n}} = {\cal I}$. 
By our choice ${\cal I}_{{\cal Y}_{2}} =\{[y_{1}, y_{2}], [y_{2},y_{1}]\}$. Hence ${\cal I}_{{\cal Y}_{j}} \neq \emptyset$ for all $j \in \{2, \ldots, n\}$. For $i \in \{2, \ldots, n-1\}$, by Proposition \ref{pro2.6} (for $\kappa = i + 1$), 
$$
{\mathcal{L}}_{\rm gr}(F_{i+1}) = {\mathcal{L}}_{\rm gr}(F_{i}) \oplus L(\{y_{i+1}\} \wr {\cal Y}_{i}^{*}). \eqno(3.1)
$$ 
As observed in \cite[Corollary II.6]{dk}, by using Lemma \ref{le2.2}, we have, for all $i \in \{2, \ldots, n-1\}$, $
I_{{\cal Y}_{i+1}} = (I_{{\cal Y}_{i+1}} \cap {\mathcal{L}}_{\rm gr}(F_{i})) \oplus (I_{{\cal Y}_{i+1}} \cap L(\{y_{i+1}\} \wr {\cal Y}_{i}^{*}))$. 
By Theorem II.7 in \cite{dk}, we get $
I_{{\cal Y}_{i+1}} \cap {\mathcal{L}}_{\rm gr}(F_{i}) = I_{{\cal Y}_{i}}$
for all $i \in \{2, \ldots, n-1\}$. Therefore 
$$
I_{{\cal Y}_{i+1}} = I_{{\cal Y}_{i}} \oplus (I_{{\cal Y}_{i+1}} \cap L(\{y_{i+1}\} \wr {\cal Y}_{i}^{*})). \eqno(3.2)
$$
For the next few lines, we recall some ideas given in \cite{dk}. Let ${\mathcal{Y}}^{*}$ be the free monoid on ${\mathcal{Y}}$; that is ${\mathcal{Y}}^{*} = \{\varepsilon\} \cup \{y_{i_{1}} \cdots y_{i_{\kappa}}: \kappa \geq 1; y_{i_{1}}, \ldots, y_{i_{\kappa}} \in {\mathcal{Y}}\}$ where $\varepsilon$ denotes the empty word. 
Write ${\mathcal{Y}}^{\kappa} = \{y_{i_{1}} \cdots y_{i_{\kappa}}: y_{i_{1}}, \ldots, y_{i_{\kappa}} \in {\mathcal{Y}}\}$ and so ${\mathcal{Y}}^{*} = \bigcup_{\kappa \geq 0}{\mathcal{Y}}^{\kappa}$ with ${\mathcal{Y}}^{0} = \{\varepsilon\}$.  
For $w = y_{i_{1}} \cdots y_{i_{s}} \in {\mathcal{Y}}^{*} \setminus \{\varepsilon\}$, the length $|w|$ of $w$ is the number $s$ of elements of ${\mathcal{Y}}$ occurring in the product of $w$. 
For a subset ${\mathcal{X}}$ of ${\mathcal{Y}}$, we denote by $|w|_{{\mathcal{X}}}$ the number of elements of ${\mathcal{X}}$ occurring in the product of $w$. 
Thus $|w| = \sum_{y \in {\mathcal{Y}}}|w|_{y}$. We define $||w||$, the multidegree of $w$, to be the $r$-tuple $(|w|_{y_{1}}, \ldots, |w|_{y_{m}})$. If $w, w^{\prime} \in {\mathcal{Y}}^{*}$, we say that $w \sim w^{\prime}$ if there exist $w_{1}, w_{2} \in {\mathcal{Y}}^{*}$ and $(a,b) \in \theta$ such that $w = w_{1}y_{a}y_{b}w_{2}$ and $w^{\prime} = w_{1}y_{b}y_{a}w_{2}$. 
We then define an equivalence relation on ${\mathcal{Y}}^{*}$ by saying that $w \equiv_{\theta}w^{\prime}$ if there exist $w_{1}, w_{2}, \ldots, w_{r} \in {\mathcal{Y}}^{*}$ such that $w = w_{1} \sim w_{2} \sim \cdots \sim w_{r-1} \sim w_{r} = w^{\prime}$. 
Let $[w]_{\theta}$ be the equivalence class of $w$ under the equivalence relation $\equiv_{\theta}$. If $w \equiv_{\theta} w^{\prime}$, then $|w| = |w^{\prime}|$ and $||w|| = ||w^{\prime}||$. 
For $w \in {\mathcal{Y}}^{*}$, let ${\mathcal{M}}_{w,\theta} = \{y \in {\mathcal{Y}}: \exists u \in {\mathcal{Y}}^{*}, [w]_{\theta} = [yu]_{\theta}\}$. 
Notice that, for $w_{1}, w_{2} \in {\mathcal{Y}}^{*}$, $w_{1} \equiv_{\theta} w_{2}$ implies ${\mathcal{M}}_{w_{1},\theta} = {\mathcal{M}}_{w_{2},\theta}$. Let $S_{\theta}$ be a system of representatives of the equivalence classes for $\equiv_{\theta}$. For $w \in {\mathcal{Y}}^{*}$, let $s(w)$ be the unique equivalent word for $\equiv_{\theta}$ to $w$ in $S_{\theta}$. 
Fix $\lambda \in \{2, \ldots, n\}$ and write $\equiv_{\theta_{\lambda}}$ for the equivalence relation induced by $\equiv_{\theta}$ on ${\mathcal{Y}}_{\lambda}^{*}$; that is, if $w, w^{\prime} \in {\mathcal{Y}}_{\lambda}^{*}$, then $w \sim w^{\prime}$ if there exist $w_{1}, w_{2} \in {\mathcal{Y}}^{*}_{\lambda}$ and $(a,b) \in \theta$ such that $w = w_{1}y_{a}y_{b}w_{2}$ and $w^{\prime} = w_{1} y_{b}y_{a}w_{2}$ with $y_{a}, y_{b} \in {\mathcal{Y}}_{\lambda}$. 
Further $w \equiv_{\theta_{\lambda}} w^{\prime}$ if there exist $w_{1}, \ldots, w_{r} \in {\mathcal{Y}}^{*}_{\lambda}$ such that $w = w_{1} \sim w_{2} \sim \cdots \sim w_{r} = w^{\prime}$. For $w \in {\mathcal{Y}}_{\lambda}^{*}$, we write $[w]_{\theta_{\lambda}}$ for the equivalence class of $w$ for $\equiv_{\theta_{\lambda}}$ and ${\mathcal{M}}_{w,\theta_{\lambda}} = \{y \in {\mathcal{Y}}_{\lambda}: \exists u \in {\mathcal{Y}}_{\lambda}^{*}, [w]_{\theta_{\lambda}} = [yu]_{\theta_{\lambda}}\}$. Further we define $S_{\theta_{\lambda}} = {\mathcal{Y}}_{\lambda}^{*} \cap S_{\theta}$. For $w \in {\mathcal{Y}}_{\lambda}^{*}$, $s_{\lambda}(w)$ denotes the unique equivalent word for $\equiv_{\theta_{\lambda}}$ to $w$ in $S_{\theta_{\lambda}}$. Thus $S_{\theta_{\lambda}} = \{s_{\lambda}(w) \in {\mathcal{Y}}_{\lambda}^{*}: s_{\lambda}(w) \in [w]_{\theta_{\lambda}}, w \in {\mathcal{Y}}_{\lambda}^{*} \}$. 
Clearly, for $\lambda = n$, $S_{\theta_{n}} = S_{\theta}$.

For $\lambda \in \{2, \ldots, n\}$ and $w = y_{i_{1}} \cdots y_{i_{r}} \in {\mathcal{Y}}_{\lambda - 1}^{r}$, we write $[y_{\lambda};w] = [y_{\lambda}, y_{i_{1}}, \ldots, y_{i_{r}}]$ and, for $w = \varepsilon$, $[y_{\lambda};\varepsilon] = y_{\lambda}$. For $\mu \geq 2$, we define 
$$
\begin{array}{rll}
{\mathcal{T}}^{\mu}_{1,\lambda} & = & \{[y_{\lambda};v]: v \in S_{\theta_{\lambda - 1}} \cap {\mathcal{Y}}^{\mu-1}_{\lambda-1}, \exists y_{r} \in {\mathcal{M}}_{v,\theta_{\lambda - 1}}, (\lambda, r) \in \theta \}, \\
{\mathcal{T}}^{\mu}_{2,\lambda} & = & \{[y_{\lambda};u]- [y_{\lambda};v]: v \in S_{\theta_{\lambda - 1}} \cap {\mathcal{Y}}^{\mu-1}_{\lambda-1}, u \in [v]_{\theta_{\lambda - 1}} \setminus \{v\}\}
\end{array}
$$
and
$$
\begin{array}{rll}
{\mathcal{T}}^{\mu}_{3,\lambda} & = & \{[y_{\lambda};v]: v \in S_{\theta_{\lambda - 1}} \cap {\mathcal{Y}}^{\mu-1}_{\lambda-1}, \forall y_{r} \in {\mathcal{M}}_{v,\theta_{\lambda - 1}}, (\lambda, r) \notin \theta \}.
\end{array}
$$
Notice that ${\mathcal{T}}^{\mu}_{1,2} = \{[y_{2}, ~_{\mu-1}y_{1}]\}$, ${\mathcal{T}}^{\mu}_{2,2} = \emptyset$ for all $\mu \geq 2$, ${\mathcal{T}}^{2}_{3,2} = \{y_{2}\}$ and ${\mathcal{T}}^{\mu}_{3,2} = \emptyset$ for all $\mu \geq 3$. Further we define ${\mathcal{T}}_{1,\lambda} = \bigcup_{\mu \geq 2}{\mathcal{T}}^{\mu}_{1,\lambda}$, ${\mathcal{T}}_{2,\lambda} = \bigcup_{\mu \geq 2}{\mathcal{T}}^{\mu}_{2,\lambda}$, ${\mathcal{T}}_{3,\lambda} = \bigcup_{\mu \geq 2}{\mathcal{T}}^{\mu}_{3,\lambda}$ and ${\mathcal{T}}^{2}_{3,\lambda} = \{y_{\lambda}\}$ for all $\lambda \geq 3$.  One can easily see that ${\cal T}_{i,\l}={\cal T}_i$ in \cite{dk}, with $i=1,2,3$.

\begin{proposition}\label{pro4.1}  
For any $\nu \in \{2, \ldots, n\}$, with $n \geq 2$, 
$$
{\mathcal{L}}_{\rm gr}(F_{\nu}) = \langle \overline{{\mathcal{F}}}_{\nu} \rangle \oplus (\bigoplus_{j=3}^{\nu}L(({\mathcal{T}}_{3,j} \setminus \{y_{j}\}) \wr \{y_{j}\}^{*}) \oplus I_{{\mathcal{Y}}_{\nu}},
$$ 
where $\langle \overline{{\mathcal{F}}}_{\nu} \rangle$ is the $\mathbb{Z}$-module spanned by $\overline{{\mathcal{F}}}_{\nu}$ and $
I_{{\mathcal{Y}}_{\nu}} =  \bigoplus_{j=2}^{\nu}L(({\mathcal{T}}_{1,j} \cup {\mathcal{T}}_{2,j}) \wr {\mathcal{T}}_{3,j}^{*})$.        
\end{proposition}

\pf We induct on $\nu$ and let $\nu = 2$. By applying Proposition \ref{pro2.6} (for $\kappa = 2$) on ${\mathcal{L}}_{\rm gr}(F_{2})$ and Theorem  \ref{th1} on  $L({\mathcal{T}}_{3,2} \bigcup {\mathcal{T}}_{1,2})$, we obtain 
$$
{\mathcal{L}}_{\rm gr}(F_{2}) = \langle y_{1} \rangle \oplus L({\mathcal{T}}_{3,2}) \oplus L({\mathcal{T}}_{1,2} \wr {\mathcal{T}}_{3,2}^{*}) = \langle y_{1}, y_{2} \rangle \oplus L({\mathcal{T}}_{1,2} \wr {\mathcal{T}}_{3,2}^{*}).
$$ 
Since ${\mathcal{T}}_{1,2} \wr {\mathcal{T}}_{3,2}^{*} \subset \gamma_{2}({\mathcal{L}}_{\rm gr}(F_{2}))$ and $\gamma_{2}({\mathcal{L}}_{\rm gr}(F_{2}))$ is an ideal of ${\mathcal{L}}_{\rm gr}(F_{2})$, we have $L({\mathcal{T}}_{1,2} \wr {\mathcal{T}}_{3,2}^{*}) \subseteq \gamma_{2}({\mathcal{L}}_{\rm gr}(F_{2}))$. 
By the modular law and since $\langle y_{1}, y_{2} \rangle \cap \gamma_{2}({\mathcal{L}}_{\rm gr}(F_{2})) = \{0\}$, we get $\gamma_{2}({\mathcal{L}}_{\rm gr}(F_{2})) = L({\mathcal{T}}_{1,2} \wr {\mathcal{T}}_{3,2}^{*})$. Since $\gamma_{2}({\mathcal{L}}_{\rm gr}(F_{2}))$ is generated as an ideal of ${\mathcal{L}}_{\rm gr}(F_{2})$ by the set ${\mathcal{I}}_{\mathcal{Y}_{2}}$, we have $I_{\mathcal{Y}_{2}} = \gamma_{2}({\mathcal{L}}_{\rm gr}(F_{2}))$ and so $I_{\mathcal{Y}_{2}} = L({\mathcal{T}}_{1,2} \wr {\mathcal{T}}_{3,2}^{*})$. 
Therefore we obtain the required result for $\nu = 2$. 

Let $\kappa \in \{2,\ldots, n - 1\}$ and assume that our claim  is valid for $\kappa$. 
Hence  ${\mathcal{L}}_{\rm gr}(F_{\kappa}) = \langle \overline{{\mathcal{F}}}_{\kappa} \rangle \oplus (\bigoplus_{j=3}^{\kappa}L({\mathcal{T}}_{3,j} \setminus \{y_{j}\}) \wr \{y_{j}\}^{*})) \oplus I_{{\mathcal{Y}}_{\kappa}}$ 
and $
I_{{\cal Y}_{\kappa}} = \bigoplus_{j = 2}^{\kappa}L(({\mathcal{T}}_{1,j} \cup {\mathcal{T}}_{2,j}) \wr {\mathcal{T}}_{3,j}^{*})$. 
By Eq. $(3.1)$ (for $i = \kappa$), ${\mathcal{L}}_{{\rm gr}}(F_{\kappa + 1}) = {\mathcal{L}}_{{\rm gr}}(F_{\kappa}) \oplus L(\{y_{\kappa + 1}\} \wr {\mathcal{Y}}_{\kappa}^{*})$. 
As shown in \cite[Theorem II.7, Corollary II.12 and Remark]{dk}, 
$$
L(\{y_{\kappa + 1}\} \wr {\mathcal{Y}}_{\kappa}^{*}) = L(\bigcup_{j=1}^{3}{\mathcal{T}}_{j,\kappa + 1}) = L({\mathcal{T}}_{3,\kappa + 1}) \oplus L(({\mathcal{T}}_{1,\kappa + 1} \cup {\mathcal{T}}_{2,\kappa + 1}) \wr {\mathcal{T}}_{3,\kappa + 1}^{*}) \eqno(3.3)
$$
and
$$
I_{{\cal Y}_{\kappa+1}} \cap L(\{y_{\kappa+1}\} \wr {\cal Y}_{\kappa}^{*}) = L(({\mathcal{T}}_{1,\kappa + 1} \cup {\mathcal{T}}_{2,\kappa + 1}) \wr ({\mathcal{T}}_{3,\kappa + 1})^{*}). \eqno(3.4)
$$

By our inductive hypothesis and Eqs $(3.3)$, $(3.4)$ and $(3.2)$, we have $
{\mathcal{L}}_{\rm gr}(F_{\kappa+1}) = \langle \overline{{\mathcal{F}}}_{\kappa+1} \rangle \oplus (\bigoplus_{j=3}^{\kappa + 1}L(({\mathcal{T}}_{3,j} \setminus \{y_{j}\}) \wr \{y_{j}\}^{*}) \oplus I_{{\mathcal{Y}}_{\kappa +1}}$. Hence, for any $\nu \in \{2, \ldots, n\}$, we obtain the desired result. \qed

\vskip .120 in

By Theorem \ref{th1}, for any $j \in \{3, \ldots, n\}$, $L({\mathcal{T}}_{3,j}) = \langle y_{j} \rangle \oplus L(({\mathcal{T}}_{3,j} \setminus \{y_{j}\}) \wr \{y_{j}\}^{*})$. By Proposition \ref{pro4.1} and since ${\mathcal{L}}_{\rm gr}(F_{n}) = {\rm gr}(F_{n})$, we obtain the following.

\begin{corollary}\label{co4.2}
For $n \geq 2$, ${\rm gr}(F_{n}) = \langle \overline{{\mathcal{F}}}_{n} \rangle \oplus (\bigoplus_{j=3}^{n}L(({\mathcal{T}}_{3,j} \setminus \{y_{j}\}) \wr \{y_{j}\}^{*}) \oplus I$, where 
$\langle \overline{{\mathcal{F}}}_{n} \rangle$ is the $\mathbb{Z}$-module spanned by 
$\overline{{\mathcal{F}}}_{n}$ and $I = \bigoplus_{j=2}^{n}L(({\mathcal{T}}_{1,j} \cup {\mathcal{T}}_{2,j}) \wr {\mathcal{T}}_{3,j}^{*})$.
\end{corollary}

Since ${\rm gr}_{d}(F_{n})$ (with $d \geq 1$) is spanned by the the set $\{[y_{i_{1}}, \ldots, y_{i_{d}}]: i_{1}, \ldots, i_{d} \in \{1, \ldots, n\}\}$, $I$ is generated by ${\mathcal{R}}_{L}$ and using the Jacobi identity in the form $[x,[y,z]] = [x,y,z] - [x,z,y]$, we have $I = \bigoplus_{d \geq 1}[R_{L}, ~_{(d-1)}{\rm gr}_{1}(F_{n})]$ where $R_{L}$ is the $\mathbb{Z}$-module spanned by ${\mathcal{R}}_{L}$. For $d \geq 2$, let $N_{d,1}$ be the subgroup of $N \cap \gamma_{d}(F_{n})$ generated by $\{[r_{a,b}, f_{i_{1}}, \ldots, f_{i_{d-2}}]: a, b, i_{1}, \ldots, i_{d-2} \in \{1, \ldots, n\}; a > b; (a,b) \in \theta\}$. Since $[R_{L}, ~_{(d-1)}{\rm gr}_{1}(F_{n})] = N_{d,1}\gamma_{d+1}(F_{n})/\gamma_{d+1}(F_{n})$, we have by Corollary \ref{co4.2}, for $d \geq 2$,
$$
I \cap {\rm gr}_{d}(F_{n}) = \bigoplus_{j=2}^{n}L^{d}_{\rm grad}(({\mathcal{T}}_{1,j} \cup {\mathcal{T}}_{2,j}) \wr {\mathcal{T}}_{3,j}^{*}) = N_{d,1}\gamma_{d+1}(F_{n})/\gamma_{d+1}(F_{n}) \eqno(3.5)
$$

\begin{proposition}\label{pro3.3}
${\mathcal{L}}(N) = I$. 
\end{proposition}

\pf Since $I = \bigoplus_{c \geq 1}(I \cap {\rm gr}_{c}(F_{n}))$ and ${\rm gr}(F_{n}) = \bigoplus_{c \geq 1}{\rm gr}_{c}(F_{n})$, we have ${\rm gr}(F_{n})/I = \bigoplus_{c \geq 1}(\frac{{\rm gr}_{c}(F_{n})+I}{I})$. Since ${\rm gr}(F_{n})/I \cong {\rm gr}(F_{n}/N)$ as Lie algebras (see \cite[Theorem 2.1]{dk2}, \cite[Theorem 6.3]{wade}), we get, for $c \geq 1$, $\frac{{\rm gr}_{c}(F_{n})+I}{I} \cong {\rm gr}_{c}(F_{n}/N)$ as $\mathbb{Z}$-modules. By Corollary \ref{co4.2}, ${\rm gr}(F_{n})/I$ is a free $\mathbb{Z}$-module and so, for $c \geq 1$, ${\rm gr}_{c}(F_{n}/N)$ is a free abelian group of finite rank. For $c \geq 1$, 
$$
\begin{array}{rll}
{\rm gr}_{c}(F_{n}/N) & \cong & (\gamma_{d}(F_{n}) \gamma_{c+1}(F_{n})N)/(\gamma_{c+1}(F_{n})N) \\
& \cong & \gamma_{c}(F_{n})/(\gamma_{c}(F_{n}) \cap (\gamma_{c+1}(F_{n})N)).
\end{array}
$$
Since $\gamma_{c+1}(F_{n}) \subseteq \gamma_{c}(F_{n})$, we have, by the modular law, 
$$
\gr_c(F_n/N)=\gamma_{c}(F_{n})/(\gamma_{c}(F_{n}) \cap (\gamma_{c+1}(F_{n})N)) = \gamma_{c}(F_{n})/(\gamma_{c+1}(F_{n})(N \cap \gamma_{c}(F_{n})) \cong {\rm gr}_{c}(F_{n})/{\rm I}_{c}(N)
$$ 
and so ${\rm gr}_{c}(F_{n}/N) \cong {\rm gr}_{c}(F_{n})/{\rm I}_{c}(N)$ for all $c$ in a natural way.  For $c \geq 1$, 
$$
\frac{{\rm gr}_{c}(F_{n})+I}{I} \cong \frac{{\rm gr}_{c}(F_{n})}{I \cap {\rm gr}_{c}(F_{n})} \cong \frac{{\rm gr}_{c}(F_{n})}{{\rm I}_{c}(N)}
$$ 
as free $\mathbb{Z}$-modules. Since $I \cap {\rm gr}_{c}(F_{n}) \subseteq {\rm I}_{c}(N) \subseteq {\rm gr}_{c}(F_{n})$ for all $c \geq 2$, we obtain $I \cap {\rm gr}_{c}(F_{n}) = {\rm I}_{c}(N)$ for all $c \geq 2$. Hence ${\mathcal{L}}(N) = I$. \qed   

%%%%%%%%%%%%%%%%%%%%%%%%%%%%%%%%%
\section{Presentation of ${\rm gr}({\rm FP(H)})$}\label{sect4}
%%%%%%%%%%%%%%%%%%%%%%%%%%%%%%%%%

%Throughout this section, we mainly use the notation of Sections \ref{sect2} and \ref{sect3}. Let $F_{n+1}$ be a free group with a free generating set ${\mathcal{F}}_{n+1} = \{f_{1}, \ldots, f_{n+1}\}$, $n \geq 2$, and let $\theta$ be a non-empty partial commutation relation on $\{1, \ldots, n+1\}$, with $(i,n+1), (n+1,i) \notin \theta$ for all $i \in \{1, \ldots, n\}$. As before, we assume that $(2,1) \in \theta$. Let $M = {\mathcal{R}}^{F_{n+1}}$ be the normal subgroup of ${\mathcal{R}}$ in $F_{n+1}$, where ${\mathcal{R}} = \{[f_{a},f_{b}], [f_{n+1}, f_{i}, f_{j}], [[f_{n+1},f_{a}], [f_{n+1},f_{b}]]: i, j \in \{1, \ldots, n\}; a > b; (a,b) \in \theta\}$. Our aim in this section is to show that ${\rm gr}(F_{n+1}/M) \cong {\rm gr}(F_{n+1})/J$ as Lie algebras in a natural way.

%%%%%%%%%%%%%%%%%%%%%%%%%%%%%%%%%%%%%
%\subsection{Auxiliary results}\label{subsec4.1}
%%%%%%%%%%%%%%%%%%%%%%%%%%%%%%%%%%%%%

\subsection{Further analysis of $J$}\label{subsec4.1} 

For each $v \in {\rm gr}(F_{n+1})$, the $\mathbb{Z}$-linear mapping ${\rm ad}v: {\rm gr}(F_{n+1}) \rightarrow {\rm gr}(F_{n+1})$ is defined by $u({\rm ad}v) = [u,v]$ for all $u \in {\rm gr}(F_{n+1})$. 
In particular, ${\rm ad}v$ is a derivation of ${\rm gr}(F_{n+1})$. The set of all derivations of ${\rm gr}(F_{n+1})$ is denoted by ${\rm Der}({\rm gr}(F_{n+1}))$ and is regarded as a Lie algebra in a natural way. 
The map ${\rm ad}: {\rm gr}(F_{n+1}) \rightarrow {\rm Der}({\rm gr}(F_{n+1}))$ sending $v$ to ${\rm ad} v$ for all $v \in {\rm gr}(F_{n+1})$ is called the adjoint representation of ${\rm gr}(F_{n+1})$. 
Let $U({\rm gr}(F_{n+1}))$ be the universal enveloping $\mathbb{Z}$-algebra of ${\rm gr}(F_{n+1})$ and we regard ${\rm gr}(F_{n+1})$ is contained in $U({\rm gr}(F_{n+1}))$. 
An element of $U({\rm gr}(F_{n+1})) \setminus \{1\}$ is a $\mathbb{Z}$-linear combination of elements of the form $u_{1} u_{2} \cdots u_{m}$ where $u_{1}, \ldots, u_{m} \in {\rm gr}(F_{n+1})$. 
For an ideal $K$ of ${\rm gr}(F_{n+1})$, we write $\gamma_{2}(K) = [K,K]$ for the $\mathbb{Z}$-submodule of $K$ spanned by all $[u,v]$, with $u, v \in K$, and $\tilde{v} = v + K$ for $v \in {\rm gr}(F_{n+1})$. Since $K$ is an ideal of ${\rm gr}(F_{n+1})$, $\gamma_{2}({\rm gr}(F_{n+1}))$ is an ideal of ${\rm gr}(F_{n+1})$. 
The $\mathbb{Z}$-module $K/\gamma_{2}(K)$ becomes a (right) ${\rm gr}(F_{n+1})/K$-module via the adjoint representation of ${\rm gr}(F_{n+1})/K$. Namely,
$$
\begin{array}{rll}
(u+\gamma_{2}(K))\tilde{v} & = & u({\rm ad}v) + \gamma_{2}(K) \\
& = & [u,v] + \gamma_{2}(K)
\end{array}
$$
for all $u \in K$ and $v \in {\rm gr}(F_{n+1})$. The $\mathbb{Z}$-module $K/\gamma_{2}(K)$ becomes a (right) $U({\rm gr}(F_{n+1})/K)$-module by defining 
$$
\begin{array}{rll}
(u+\gamma_{2}(K)) \tilde{v}_{1} \tilde{v}_{2} \cdots \tilde{v}_{\kappa} & = & u({\rm ad}v_{1}) ({\rm ad}v_{2}) \cdots ({\rm ad}v_{\kappa}) + \gamma_{2}(K) \\
& = & [u, v_{1}, v_{2}, \ldots, v_{\kappa}] + \gamma_{2}(K)
\end{array}
$$
for all $u \in K$ and $v_{1}, v_{2}, \ldots, v_{\kappa} \in {\rm gr}(F_{n+1})$.

\begin{lemma}\label{le5.1}
The $\mathbb{Z}$-module $L({\cal E} \wr \Omega^{*})/\gamma_{2}(L({\cal E} \wr \Omega^{*}))$ is a free (right) $U({\rm gr}(F_{n+1})/L({\cal E} \wr \Omega^{*}))$-module with a free generating set $\{[\t, y_{i}, y_{j}] + \gamma_{2}(L({\cal E} \wr \Omega^{*})): i, j \in \n\}$. 
\end{lemma}

\pf Since $
{\cal E} \wr \Omega^{*} = \{[\t, y_{i_{1}}, \ldots, y_{i_{\kappa}}, b_{j_{1}}, \ldots, b_{j_{\lambda}}]: \kappa \geq 2; \lambda \geq 0; y_{i_{1}} \cdots y_{i_{\kappa}} \in (\overline{{\mathcal F}}_{n})^{*}; b_{j_{1}} \cdots b_{j_{\lambda}} \in \Omega^{*}\}$, 
we have the $\mathbb{Z}$-module $L({\cal E} \wr \Omega^{*})/\gamma_{2}(L({\cal E} \wr \Omega^{*}))$ is free on ${\cal E} \wr \Omega^{*}$ modulo $\gamma_{2}(L({\cal E} \wr \Omega^{*}))$. 
For the next few lines, let $e_{ij} = [\t, y_{i}, y_{j}]$ for all $i, j \in \n$ and $\widehat{J} = L({\cal E} \wr \Omega^{*})$. 
We point out that $e_{ij} = [\omega_{i},y_{j}]$ with $\omega_{i} = [\t,y_{i}]$ for all $i, j \in \n$. 
Let $
{\mathcal{L}}_{{\rm gr},\widehat{J}}(F_{n}) = \frac{\mathcal{L}_{{\rm gr}}(F_{n}) + \widehat{J}}{\widehat{J}}$ and $L_{\widehat{J}}(\Omega) = \frac{L(\Omega) + \widehat{J}}{\widehat{J}}$. 
By Eq. $(2.8)$, we obtain $
{\rm gr}(F_{n+1})/\widehat{J} = {\mathcal{L}}_{{\rm gr},\widehat{J}}(F_{n}) \oplus L_{\widehat{J}}(\Omega)$. 
Thus the free $\mathbb{Z}$-module ${\rm gr}(F_{n+1})/\widehat{J}$ is a direct sum, as $\mathbb{Z}$-modules, of its Lie subalgebras ${\mathcal{L}}_{{\rm gr},\widehat{J}}(F_{n})$ and $L_{\widehat{J}}(\Omega)$. 
We point out that ${\mathcal{L}}_{{\rm gr},\widehat{J}}(F_{n}) \cong \mathcal{L}_{{\rm gr}}(F_{n})$ and $L_{\widehat{J}}(\Omega) \cong L(\Omega)$ as Lie algebras in a natural way. 
For $i \in \{1,\ldots,r\}$, let $\overline{w}_i = w_i + \widehat{J} \in U({\rm gr}(F_{n+1})/\widehat{J}) \setminus \{0\}$ and $w_i \in {\rm gr}(F_{n+1})$. 
Write each $w_{i} = u_{i} + v_{i} + t_{i}$ with $u_{i} \in \mathcal{L}_{{\rm gr}}(F_{n})$, $v_{i} \in L(\Omega)$ and $t_{i} \in \widehat{J}$. 
Furthermore we assume that each $u_{i}$, $v_{i}$ is a simple Lie commutator of degree $\geq 1$ in terms of the elements of the sets $\overline{\mathcal{F}}_{n}$ and $\Omega$, respectively. Then, for all $i, j \in \n$,  
$$
\begin{array}{rll}
(e_{ij} + \gamma_{2}(\widehat{J})) \overline{w}_{1} \cdots \overline{w}_{r} & = & e_{ij}({\rm ad}w_{1}) \cdots ({\rm ad}w_{r}) + \gamma_{2}(\widehat{J}) \\
& = & [e_{ij}, w_{1}, \ldots, w_{r}] + \gamma_{2}(\widehat{J}) \\
& = & [e_{ij}, u_{1} + v_{1} + t_{1}, \ldots, u_{r} + v_{r} + t_{r}] + \gamma_{2}(\widehat{J}).  
\end{array}
$$
By the multi-linearity of the Lie commutator, $e_{ij}, t_{1}, \ldots, t_{r} \in \widehat{J}$ and by using the Jacobi identity in the form $[a,[b,c]] = [a,b,c] - [a,c,b]$, we can write the element $(e_{ij}+ \gamma_{2}(\widehat{J})) \overline{w}_{1} \cdots \overline{w}_{r}$ as a $\mathbb{Z}$-linear combination of elements of the form $[e_{ij}, z_{1}, \ldots, z_{\mu}] + \gamma_{2}(\widehat{J})$, where $z_{1}, \ldots, z_{\mu} \in \overline{\mathcal{F}}_{n} \cup \Omega$. 
For a moment, we consider the product $z_{1} \cdots z_{\mu}$ as a $\mu$-tuple $(z_{1}, \ldots, z_{\mu})$ and assume that, for some $j \in \n$, $y_{j}$ occurs in the $k$-th position ($k \geq 2$) of $(z_{1}, \ldots, z_{\mu})$. If in the $(k-1)$-th position of $(z_{1}, \ldots, z_{\mu})$ is $\omega_{\tau}$, for some $\tau \in \n$, then 
$$
[e_{ij},z_{1}, \ldots, \omega_{\tau},y_{j}, \ldots, z_{\mu}] + \gamma_{2}(\widehat{J}) = [e_{ij},z_{1}, \ldots, y_{j}, \omega_{\tau}, \ldots, z_{\mu}] + \gamma_{2}(\widehat{J}), \eqno(4.1)
$$
since both $e_{ij}$ and $[\omega_{\tau},y_{j}] = e_{\tau j} \in \widehat{J}$. If in the $(k-1)$-position of $(z_{1}, \ldots, z_{\mu})$ there is a $\t$, then 
$$
[e_{ij},z_{1}, \ldots, \t,y_{j}, \ldots, z_{\mu}] -  ([e_{ij},z_{1}, \ldots, y_{j}, \t, \ldots, z_{\mu}] + [e_{ij},z_{1}, \ldots, \omega_{j}, \ldots, z_{\mu}]) \in \gamma_{2}(\widehat{J}). \eqno(4.2)
$$
Having in mind $(4.1)$ and $(4.2)$, we may write each element $[e_{ij}, z_{1}, \ldots, z_{\mu}] + \gamma_{2}(\widehat{J})$, where $z_{1}, \ldots, z_{\mu} \in \overline{\mathcal{F}}_{n} \cup \Omega$, as  a $\mathbb{Z}$-linear combination of elements of the form         
$$
e_{ij}({\rm ad}y_{i_{1}}) \cdots ({\rm ad}y_{i_{\kappa}})({\rm ad}b_{j_{1}}) \cdots ({\rm ad} b_{j_{\lambda}}) + \gamma_{2}(\widehat{J}) =  (e_{ij} + \gamma_{2}(\widehat{J})) \overline{y}_{i_{1}} \cdots  \overline{y}_{i_{\kappa}} \overline{b}_{j_{1}} \cdots \overline{b}_{j_{\lambda}},
$$
where $y_{i_{1}} \cdots y_{i_{\kappa}} \in (\overline{\mathcal{F}}_{n})^{*}$ and $ b_{j_{1}} \cdots b_{j_{\lambda}} \in \Omega^{*}$. Hence, since the set ${\cal E} \wr \Omega^{*}$ is a free generating set for $\widehat{J}$, we have the required result. \qed     

\begin{proposition}\label{pro4.2} 
Let $M_{{\cal R}_{3}} = {\cal R}_{3}^{F_{n+1}}$ be the normal closure of ${\cal R}_{3}$ in $F_{n+1}$ where ${\cal R}_{3} = \{[t,f_{i},f_{j}]: i, j \in \{1, \ldots, n\}\}$. Then ${\cal L}(M_{{\cal R}_{3}}) = L({\cal E} \wr \Omega^{*})$. 
\end{proposition}

\pf By Proposition \ref{pro3.6}~(4), $L({\cal E} \wr \Omega^{*})$ is the ideal of ${\rm gr}(F_{n+1})$ generated by ${\mathcal{R}}_{3,L} = \{[\t,y_{i},y_{j}]: i, j \in \{1, \ldots, n\}\}$. 
Since ${\cal R}_{3,L} \subset {\rm I}_{3}(M_{{\cal R}_{3}})$ and ${\cal L}(M_{{\cal R}_{3}})$ is an ideal of ${\rm gr}(F_{n+1})$, we have $L({\cal E} \wr \Omega^{*}) \subseteq {\cal L}(M_{{\cal R}_{3}})$. 
By Eq. $(2.8)$, $L({\cal E} \wr \Omega^{*})$ is a direct summand of ${\rm gr}(F_{n+1})$ and so ${\rm gr}(F_{n+1})/L({\cal E} \wr \Omega^{*})$ is a free $\mathbb{Z}$-module. 
By Lemma \ref{le5.1}, $L({\cal E} \wr \Omega^{*})/\gamma_{2}(L({\cal E} \wr \Omega^{*}))$ is a free (right) $U({\rm gr}(F_{n+1})/L({\cal E} \wr \Omega^{*}))$-module with a free generating set ${\cal R}_{3,L}+\gamma_{2}(L({\cal E} \wr \Omega^{*}))$. 
By a result of Labute \cite[Theorem 1]{labu}, $
{\rm gr}(F_{n+1}/M_{{\cal R}_{3}}) \cong {\rm gr}(F_{n+1})/L({\cal E} \wr \Omega^{*})$   
as Lie algebras. 
Since $L({\cal E} \wr \Omega^{*}) = \bigoplus_{c \geq 1}(L({\cal E} \wr \Omega^{*}) \cap {\rm gr}_{c}(F_{n+1}))$ and ${\rm gr}(F_{n+1}) = \bigoplus_{c \geq 1}{\rm gr}_{c}(F_{n+1})$, we have ${\rm gr}(F_{n+1})/L({\cal E} \wr \Omega^{*}) = \bigoplus_{c \geq 1}(\frac{{\rm gr}_{c}(F_{n+1})+L({\cal E} \wr \Omega^{*})}{L({\cal E} \wr \Omega^{*})})$. 
Since ${\rm gr}(F_{n+1})/L({\cal E} \wr \Omega^{*}) \cong {\rm gr}(F_{n+1}/M_{{\mathcal{R}_{3}}})$, we get, for $c \geq 1$, $\frac{{\rm gr}_{c}(F_{n+1})+L({\cal E} \wr \Omega^{*})}{L({\cal E} \wr \Omega^{*})} \cong {\rm gr}_{c}(F_{n+1}/M_{\mathcal{R}_{3}})$ as $\mathbb{Z}$-modules. 
By Eq. $(2.8)$, ${\rm gr}(F_{n+1})/L({\cal E} \wr \Omega^{*})$ is a free $\mathbb{Z}$-module and so, for $c \geq 1$, ${\rm gr}_{c}(F_{n+1}/M_{\mathcal{R}_{3}})$ is a free abelian group of finite rank. 
As in the proof of Proposition \ref{pro3.3}, for $c \geq 1$, $
{\rm gr}_{c}(F_{n+1}/M_{\mathcal{R}_{3}}) \cong {\rm gr}_{c}(F_{n+1})/{\rm I}_{c}(M_{\mathcal{R}_{3}})$ in a natural way. 
For $c \geq 1$, 
$$
\frac{{\rm gr}_{c}(F_{n+1})+L({\cal E} \wr \Omega^{*})}{L({\cal E} \wr \Omega^{*})} \cong \frac{{\rm gr}_{c}(F_{n+1})}{{\rm gr}_{c}(F_{n+1}) \cap L({\cal E} \wr \Omega^{*})} \cong \frac{{\rm gr}_{c}(F_{n+1})}{{\rm I}_{c}(M_{\mathcal{R}_{3}})}
$$ 
as free $\mathbb{Z}$-modules. Since ${\rm gr}_{c}(F_{n+1}) \cap L({\cal E} \wr \Omega^{*}) \subseteq {\rm I}_{c}(M_{\mathcal{R}_{3}}) \subseteq {\rm gr}_{c}(F_{n+1})$ for all $c \geq 2$, we obtain ${\rm gr}_{c}(F_{n+1}) \cap L({\cal E} \wr \Omega^{*}) = {\rm I}_{c}(M_{\mathcal{R}_{3}})$ for all $c \geq 2$. 
Hence we obtain  ${\cal L}(M_{{\cal R}_{3}}) = L({\cal E} \wr \Omega^{*})$. \qed

\begin{proposition}\label{pro4.3}
Let $M_{\mathcal{R}_{2}} = {\mathcal{R}}_{2}^{F_{n}}$ be the normal closure of ${\mathcal{R}}_{2}$ in $F_{n}=\<f_1,\ldots,f_n\>$ where ${\mathcal{R}}_{2} = \{[f_{a},f_{b}]: a > b; (a,b) \in \theta\}$. Then ${\mathcal{L}}(M_{\mathcal{R}_{2}}) = J \cap {\mathcal{L}}_{\rm gr}(F_{n})$.

\end{proposition}

\pf Recall that $N_{t}$ is the normal closure of $\{t\}$ in $F_{n+1}$. Since $N_{t} \cap M_{\mathcal{R}_{2}} = \{1\}$, we have $M_{\mathcal{R}_{2}} \cap \gamma_{d}(F_{n+1}) = M_{\mathcal{R}_{2}} \cap \gamma_{d}(F_{n})$ for all $d \geq 1$. Hence, for $d \geq 1$, ${\rm I}_{d}(M_{\mathcal{R}_{2}}) = \frac{(M_{\mathcal{R}_{2}} \cap \gamma_{d}(F_{n})) \gamma_{d+1}(F_{n+1})}{\gamma_{d+1}(F_{n+1})} \leq {\mathcal{L}}_{{\rm gr},d}(F_{n})$. Since $M_{\mathcal{R}_{2}}$ is normal in $F_{n}$, we obtain ${\mathcal{L}}(M_{\mathcal{R}_{2}}) = \bigoplus_{d \geq 2}{\rm I}_{d}(M_{\mathcal{R}_{2}})$ is an ideal in ${\mathcal{L}}_{\rm gr}(F_{n})$. By Proposition \ref{pro3.6}~(2), $J \cap {\mathcal{L}}_{\rm gr}(F_{n}) \subseteq {\mathcal{L}}(M_{\mathcal{R}_{2}})$. By Proposition \ref{pro3.3} (for $N = M_{\mathcal{R}_{2}}$, $I = J \cap {\mathcal{L}}_{\rm gr}(F_{n})$ and ${\mathcal{L}}_{\rm gr}(F_{n}) = {\rm gr}(F_{n})$), we get ${\mathcal{L}}(M_{\mathcal{R}_{2}}) = J \cap {\mathcal{L}}_{\rm gr}(F_{n})$. \qed

\bigskip

For $i \in \n$, let $q_{i} = [t,f_{i}]$,  $\omega_{i} = q_{i}\gamma_{3}(F_{n+1}) = [\t,y_{i}]$ and  $\Omega_{0} = \{\omega_{1}, \ldots, \omega_{n}\}$. 
Thus $\Omega_{0} \subset {\rm gr}_{2}(F_{n+1})$ and $\Omega = \{\t\} \cup \Omega_{0}$. By Theorem \ref{th1}, $L(\Omega) = L(\Omega_{0}) \oplus L(\{\t\} \wr \Omega^{*}_{0})$. 
We write $F_{Q}$ for the subgroup of $F_{n+1}$ generated by $Q = \{q_{1}, \ldots, q_{n}, t\}$. 
It is clear enough that $F_{Q}$ is a free group of rank $n+1$. In particular $F_{n+1}$ is obviously isomorphic to $F_{Q}$. Let $F_{Q_{0}}$ be the  subgroup of $F_{Q}$ generated by $Q_{0} = \{q_{1}, \ldots, q_{n}\}$. Clearly $F_{Q_{0}}$ is free on $Q_{0}$. 
We point out that $F_{Q}$ is the semidirect product of $N_{Q,t}$ by $F_{Q_{0}}$ where $N_{Q,t}$ is the normal closure of $\{t\}$ in $F_{Q}$. 
For $c \geq 1$, let $
{\cal L}_{{\rm gr},c}(F_{Q_{0}}) = \gamma_{c}(F_{Q_{0}})\gamma_{2c+1}(F_{n+1})/\gamma_{2c+1}(F_{n+1})$  
and form the (restricted) direct sum ${\cal L}_{{\rm gr}}(F_{Q_{0}})$ of the abelian groups ${\cal L}_{{\rm gr},c}(F_{Q_{0}})$. 
Since $[\gamma_{\k}(F_{Q_{0}}), \gamma_{\l}(F_{Q_{0}})] \subseteq \gamma_{\k+\l}(F_{Q_{0}})$ for all $\k,\l$, we have ${\cal L}_{\rm gr}(F_{Q_{0}})$ is a Lie subalgebra of ${\rm gr}(F_{n+1})$. 
Notice that, for $c \geq 1$, ${\cal L}_{{\rm gr},c}(F_{Q_{0}}) = L^{c}({\Omega_{0}})$ and so ${\cal L}_{\rm gr}(F_{Q_{0}}) = L({\Omega_{0}})$. 
Therefore $L(\Omega) = {\cal L}_{\rm gr}(F_{Q_{0}}) \oplus L(\{\t\} \wr \Omega^{*}_{0})$. For $c \geq 1$, let $F_{Q,c} = F_{Q} \cap \gamma_{c}(F_{n+1})$. 
Form the (restricted) direct sum ${\cal L}(F_{Q})$ of the abelian groups ${\rm I}_{c}(F_{Q}) = F_{{Q},c}\gamma_{c+1}(F_{n+1})/\gamma_{c+1}(F_{n+1})$. 
Clearly ${\cal L}(F_{Q})$ is a Lie subalgebra of ${\rm gr}(F_{n+1})$. 

\begin{lemma}\label{le4.4}  ${\cal L}(F_{Q}) = L(\Omega)$.
\end{lemma}

\pf Since $F_{Q}$ is the semidirect product of $N_{Q,t}$ by $F_{Q_{0}}$ and $F_{{Q_{0}}} \subseteq F^{\prime}_{n+1}$, we have ${\rm I}_{1}(F_{Q}) = {\rm I}_{1}(N_t)$ and ${\rm I}_{2}(F_{Q}) = {\rm I}_{2}(N_t) = {\mathcal{L}}_{{\rm gr},1}(F_{Q_{0}})$. 
Notice that $L(\Omega) = L({\rm I}_{1}(N_t) \oplus {\rm I}_{2}(N_t))$. Since ${\cal L}(F_{Q})$ is a Lie subalgebra of ${\rm gr}(F_{n+1})$ and $L(\Omega)$ is free on ${\rm I}_{1}(N_t) \oplus {\rm I}_{2}(N_t)$, we have $L(\Omega) \subseteq {\cal L}(F_{Q})$. Clearly ${\cal L}(F_{Q}) \subseteq {\cal L}(N_t)$. 
By Eq. $(2.6)$ and the modular law, we get 
$$
{\cal L}(F_{Q}) = L(\Omega) \oplus (L({\cal E} \wr \Omega^{*}) \cap {\cal L}(F_{Q})). \eqno (4.3)
$$
We claim that $L({\cal E} \wr \Omega^{*}) \cap {\cal L}(F_{Q}) = \{0\}$. 
Let $\overline{u} \in L({\cal E} \wr \Omega^{*}) \cap {\cal L}(F_{Q})$. By Proposition \ref{pro4.2}, ${\cal L}(M_{{\cal R}_{3}}) = L({\cal E} \wr \Omega^{*})$. 
Since both ${\cal L}(M_{{\cal R}_{3}})$ and ${\cal L}(F_{Q})$ are graded Lie algebras, we may assume that $\overline{u} = u \gamma_{c+1}(F_{n+1}) \in {\rm I}_{c}(M_{{\cal R}_{3}}) \cap {\rm I}_{c}(F_{Q})$ with $u \in \gamma_{c}(F_{n+1})$ and some $c \geq 1$. 
It is clear enough that $c \geq 3$. 
Then $\overline{u} = u_{1}\gamma_{c+1}(F_{n+1}) = u_{2}\gamma_{c+1}(F_{n+1})$, 
where $u_{1} \in M_{{\mathcal{R}}_{3}} \cap \gamma_{c}(F_{n+1})$ and $u_{2} \in F_{Q,c} = F_{Q}\cap \g_c(F_{n+1})$. 
To get a contradiction, we assume that $u \notin \gamma_{c+1}(F_{n+1})$. Hence $u_{1}, u_{2} \in \gamma_{c}(F_{n+1}) \setminus \gamma_{c+1}(F_{n+1})$.  
Since $F_{Q}$ is a free group on $Q = \{q_{1}, \ldots, q_{n}, t\}$, $u_{2}$ is (uniquely) written as $u_{2} = z^{\mu_{1}}_{i_{1}} \cdots z^{\mu_{\kappa}}_{i_{\kappa}}$, with $z_{i_{1}}, \ldots, z_{i_{\kappa}} \in Q$, $i_{\lambda} \neq i_{\lambda + 1}$, $\lambda \in \{1,\ldots,\kappa - 1\}$ and $\mu_{1}, \ldots, \mu_{\kappa} \in {\mathbb{Z}} \setminus \{0\}$. 
Notice that the elements of $Q$ are ordered as $t < q_{1} < \cdots < q_{n}$. In the next few lines, we write $\widehat{f}_{i} = f_{i}\gamma_{c+1}(F_{n+1})$, $i \in \{1, \ldots, n\}$ and $\widehat{t}=t\g_{c+1}(F_{n+1})$. 
Thus $\{\widehat{f}_{1}, \ldots, \widehat{f}_{n},\widehat{t}\}$ is a free generating set of $F_{n+1,c} = F_{n+1}/\gamma_{c+1}(F_{n+1})$. 
Observe that, for any $j \in \n$, $\widehat{q}_{j} = q_{j}\gamma_{c+1}(F_{n+1}) = [\widehat{t}, \widehat{f}_{j}] \in \gamma_{2}(F_{n+1,c})$. 
Using the collection process for arbitrary group elements on $u_{2} = z^{\mu_{1}}_{i_{1}} \cdots z^{\mu_{\kappa}}_{i_{\kappa}}$ (see \cite[Chapter 3, Theorem 3.1 and Theorem 3.5]{cmz}), the element $\widehat{u}_{2} = u_{2}\gamma_{c+1}(F_{n+1})$ is written in $F_{n+1,c}$ as 
$$
\widehat{u}_{2} = (\widehat{t})^{\alpha}(\widehat{q}_{1})^{\alpha_{1}} \cdots (\widehat{q}_{n})^{\alpha_{n}} w_{3}(\widehat{t}, \widehat{q}_{1}, \ldots, \widehat{q}_{n}) \cdots w_{c}(\widehat{t}, \widehat{q}_{1}, \ldots, \widehat{q}_{n}),
$$
where $\alpha, \alpha_{1}, \ldots, \alpha_{n}$ are integers, $w_{3}(t,q_{1}, \ldots, q_{n}), \ldots, w_{c}(t,q_{1}, \ldots,q_{n})$ are products of basic group commutators in $t$, $f_{1}, \ldots, f_{n}$ and have weights $3, 4, \ldots, c$ in $F_{n+1}$, respectively.  
Thus, for $j \in \{3,\ldots,c\}$, $w_{j}(t, q_{1}, \ldots, q_{n}) \in F_{Q,j}$. 
Since $u_{2} \in \gamma_{c}(F_{n+1}) \setminus \gamma_{c+1}(F_{n+1})$ and the center of $F_{n+1,c}$ equals $\gamma_{c}(F_{n+1,c})$,  we get $\widehat{u}_{2} = w_{c}(\widehat{t}, \widehat{q}_{1}, \ldots, \widehat{q}_{n})$. Therefore $u_{2}\gamma_{c+1}(F_{n+1}) \in L^{c}_{\rm grad}(\Omega)$. 
On the other hand, $u_{2}\gamma_{c+1}(F_{n+1}) = u_{1}\gamma_{c+1}(F_{n+1}) \in {\rm I}_{c}(M_{{\cal R}_{3}}) = L^{c}_{\rm grad}(\mathcal{E} \wr \Omega^{*})$ which contradicts to the fact $L(\Omega) \cap L(\mathcal{E} \wr \Omega^{*}) = \{0\}$ (Eq. (2.8)). 
Therefore $L({\cal E} \wr \Omega^{*}) \cap {\cal L}(F_{Q}) = \{0\}$ and so, by Eq. $(4.3)$, ${\cal L}(F_{Q}) = L(\Omega)$. \qed        
            
\begin{lemma}\label{le4.5}
$J \cap L(\Omega)$ is a direct summand of ${\cal L}(F_{Q})$.
\end{lemma}

\pf Let $\psi$ be the Lie algebra isomorphism from ${\rm gr}(F_{n+1})$ onto ${\mathcal{L}}(F_{Q})$ satisfying the conditions $\psi(\t) = \t$ and $\psi(y_{i}) = \omega_{i} = [\t,y_{i}]$, $i \in \{1, \ldots, n\}$. Let $J_{1}$ be the ideal of ${\rm gr}(F_{n+1})$ generated by ${\mathcal{J}}_{1} = \{[y_{a},y_{b}]: a > b; (a,b) \in \theta\}$.  Since $\psi$ is a Lie algebra isomorphism and Proposition \ref{pro3.6}~(3), we have $\psi(J_{1}) = J \cap L(\Omega)$. By Corollary \ref{co4.2} (for $n+1$), we get
$$
{\rm gr}(F_{n+1}) = \langle \overline{\mathcal{F}}_{n+1}\rangle \oplus (\bigoplus_{j=3}^{n+1}L((\mathcal{T}_{3,j} \setminus \{y_{j}\}) \wr \{y_{j}\}^{*})) \oplus J_{1},
$$
where $\langle \overline{\mathcal{F}}_{n+1}\rangle$ is the $\mathbb{Z}$-module spanned by $\overline{\mathcal{F}}_{n+1}$ and $J_{1} = \bigoplus_{j=2}^{n}L((\mathcal{T}_{1,j} \cup (\mathcal{T}_{2,j}) \wr (\mathcal{T}_{3,j}^{*}) \oplus L(\mathcal{T}_{2,n+1} \wr (\mathcal{T}_{3,n+1}^{*})$. Notice that since $(n+1,i), (i,n+1) \notin \theta$ for $i \in \{1, \ldots, n\}$, we have $\mathcal{T}_{1,n+1} = \emptyset$. 
Since $\psi$ is a Lie algebra isomorphism, we obtain 
$$
{\mathcal{L}}(F_{Q}) = \langle \overline{\mathcal{Q}} \rangle \oplus (\bigoplus_{j=3}^{n+1}L((\psi(\mathcal{T}_{3,j}) \setminus \{\omega_{j}\}) \wr \{\omega_{j}\}^{*})) \oplus \psi(J_{1}),
$$
where $\langle \overline{\mathcal{Q}}\rangle$ is the $\mathbb{Z}$-module spanned by $\overline{\mathcal{Q}} = \{\t, \omega_{1}, \ldots, \omega_{n}\}$. 
Since $\psi(J_{1}) = J \cap L(\Omega)$, we obtain the desired result. \qed      

\begin{proposition}\label{pro4.6}
Let $M_{{\cal R}_{4}} = {\cal R}_{4}^{F_{Q}}$ be the normal closure of ${\cal R}_{4}$ in $F_{Q}$, where ${\mathcal{R}}_{4} = \{[[t,f_{a}],[t,f_{b}]]: a > b; (a,b) \in \theta\}$. Then  ${\cal L}(M_{{\cal R}_{4}}) = J \cap L(\Omega)$.  
\end{proposition}   

\pf For $c \geq 1$, let $M_{{\cal R}_{4},c} = M_{{\cal R}_{4}} \cap \gamma_{c}(F_{n+1})$. 
Clearly $M_{{\cal R}_{4},c} = M_{{\cal R}_{4}}$ for $c \in \{1,2,3,4\}$. For $c \geq 4$, let ${\rm I}_{c}(M_{{\cal R}_{4}}) = (M_{{\cal R}_{4},c} \gamma_{c+1}(F_{n+1}))/\gamma_{c+1}(F_{n+1})$ and ${\cal L}(M_{{\cal R}_{4}}) = \bigoplus_{c \geq 4}{\rm I}_{c}(M_{{\cal R}_{4}})$. 
Since $M_{{\cal R}_{4},c}  \subseteq F_{{Q},c}$ for all $c$, we have ${\cal L}(M_{{\cal R}_{4}})$ is a Lie subalgebra of ${\cal L}(F_{Q})$. 
By Lemma \ref{le4.4}, ${\cal{L}}(F_{Q}) = L(\Omega)$. Since $M_{{\cal R}_{4}}$ is a normal subgroup of $F_{Q}$, we obtain ${\cal L}(M_{{\cal R}_{4}})$ is an ideal of ${\cal L}(F_{Q})$. 
By Proposition \ref{pro3.6}~(3), $J \cap L(\Omega)$ is the ideal of ${\cal L}(F_{Q})$ generated by ${\cal R}_{4,L} = \{[\omega_{a},\omega_{b}]: a > b; (a,b) \in \theta\}$ and so $J \cap L(\Omega) = \bigoplus_{d \geq 4} (J \cap {\cal L}_{{\rm gr},c}(F_{Q_{0}}))$. 
By Lemma \ref{le4.5}, ${\mathcal{L}}(F_{Q})/(J \cap L(\Omega))$ is a free $\mathbb{Z}$-module. 
By the proof of Proposition \ref{pro3.3} (for ${\rm gr}(F_{n}) = {\mathcal{L}}(F_{Q})$, $N = M_{{\cal{R}}_{4}}$ and $I = J \cap L(\Omega)$), we obtain  ${\cal L}(M_{{\cal R}_{4}}) = J \cap L(\Omega)$. \qed    

\subsection{Analysis of $M$}

The proof  of the following result is elementary and is based in calculations using properties of commutator identities.

\begin{lemma}\label{le4.7}
Let $G$ be a group and $Y$ be a subgroup of $G$. 
\begin{enumerate}

\item (Jacobi identity) For any $a \in Y$ and $b, c \in G$, $[a,b,c] = [a,c,b] [a,[b,c]] v$, where $v \in [Y,G,G,G]$.

\item For a positive integer $m$, with $m \geq 2$, and $\kappa \in \{0, \ldots, m-1\}$, let $S_{m,\kappa}$ be the the set of all permutations of $\{1, \ldots, m\}$ satisfying the conditions $\sigma(1) > \cdots > \sigma(\kappa) > \sigma(\kappa + 1) < \sigma(\kappa + 2) < \cdots < \sigma(m)$ and $S^{*}_{m} = \bigcup^{m-1}_{\kappa = 0} S_{m,\kappa}$. If $a \in Y$, $g_{1}, \ldots, g_{m} \in G$, then  
$$
[a,[g_{1}, \ldots, g_{m}]]) =\left(\prod_{\kappa = 0}^{m-1}\left(\prod_{\sigma \in S_{m,\kappa}}[a,g_{\sigma(1)}, \ldots, g_{\sigma(m)}]^{(-1)^{\kappa}}\right)\right)w,
$$
where $w \in [Y,~\prescript{}{m+1}G]$.  
\end{enumerate} 
\end{lemma}    

%Let $F_{n+1}$ be the free group of rank $n+1$, $n \geq 2$, with a free generating set ${\cal{F}}_{n+1} = \{f_{1}, \ldots, f_{n}, f_{n+1}\}$. 
Recall that ${\cal R}_{2} = \{[f_{a},f_{b}]: a > b; (a, b) \in \theta\}$,
${\cal R}_{3} = \{[t,f_{i},f_{j}]: i, j \in \n\}$,  
${\cal R}_{4} = \{[[t,f_{a}],[t,f_{b}]]: a > b; (a,b) \in \theta\}$ and ${\cal R} = {\cal R}_{2} \cup {\cal R}_{3} \cup {\cal R}_{4}$. For $j \in \{2,3,4\}$, let $R_{j}$ be the subgroup of $F_{n+1}$ generated by ${\cal R}_{j}$ and $R$ be the subgroup of $F_{n+1}$ generated by ${\cal R}$. Since $R \subseteq F^{\prime}_{n+1}$, we have $M \subseteq F^{\prime}_{n+1}$. 

\begin{lemma}\label{le4.8}
\begin{enumerate}
\item $M = R ~[R_{2}, F_{n+1}]~ [R_{3}, F_{n+1}]~ [R_{4}, F_{n+1}]$.   

\item Let $v = [r, z^{\varepsilon_{1}}_{1}, \ldots, z^{\varepsilon_{\kappa}}_{\kappa}]$ with $r \in {\cal R}_{2}$, $\kappa \geq 1$, $\varepsilon_{i} = \pm 1$, $i \in \{1,\ldots,\kappa\}$ and  $z_{1}, \ldots, z_{\kappa} \in {\cal F}_{n+1}$. 
If at least one of $z_{1}, \ldots, z_{\kappa}$ is $t$, then $v\gamma_{\kappa + 3}(F_{n+1}) \in {\rm I}_{\kappa + 2}(M_{{\cal R}_{3}})$. Otherwise $v\gamma_{\kappa + 3}(F_{n+1}) \in {\rm I}_{\kappa + 2}(M_{{\cal R}_{2}})$.   

\item Let $v = [r, z^{\varepsilon_{1}}_{1}, \ldots, z^{\varepsilon_{\kappa}}_{\kappa}]$ with $r \in {\cal R}_{4}$, $\kappa \geq 1$, $\varepsilon_{i} = \pm 1$, $i \in \{1,\ldots,\kappa\}$ and  $z_{1}, \ldots, z_{\kappa} \in {\cal F}_{n+1}$. Then $v\gamma_{\kappa + 5}(F_{n+1}) \in {\rm I}_{\kappa + 4}(M_{{\cal R}_{3}}) \oplus {\rm I}_{\kappa + 4}(M_{{\cal R}_{4}})$. 

\end{enumerate}

\end{lemma}

\pf 
\begin{enumerate}

\item Throughout this proof, we write $\widetilde{R} = R [R_{2}, F_{n+1}] [R_{3}, F_{n+1}] [R_{4}, F_{n+1}]$. Since each $[R_{j}, F_{n+1}]$ is normal in $F_{n+1}$, we have each $[v, h_{1}, \ldots, h_{\lambda}] \in [R_{j}, F_{n+1}]$ for all $v \in {\cal{R}}_{j}$ and $h_{1}, \ldots, h_{\lambda} \in F_{n+1}$. Since $M = R~[R,F_{n+1}]$,  it is enough to show that $[R, F_{n+1}] \subseteq \widetilde{R}$. In particular it is enough to show that any $[u,f] \in \widetilde{R}$ for $u \in R$ and $f \in F_{n+1} \setminus \{1\}$. Write $u = u_{i_{1}}^{\varepsilon_{1}} \cdots u_{i_{\kappa}}^{\varepsilon_{\kappa}}$ with $\kappa \geq 1$, $\varepsilon_{j} = \pm 1$, $j \in \{1, \ldots, \kappa\}$ and $u_{i_{1}}, \ldots, u_{i_{\kappa}} \in {\cal{R}}$. Using an inductive argument on $\kappa$ and the following identities 
\begin{enumerate}
\item $[ab,c] = [a,c][a,c,b][b,c]$,

\item $[a,bc] = [a,c][a,b][a,b,c]$,

\item $[a^{-1},b] = [a,b,a][a,b,a,a^{-1}][a,b]^{-1}$,

\item $[a,b^{-1}] = [a,b,b][a,b,b,b^{-1}][a,b]^{-1}$

\end{enumerate}
we write $[u,f]$ as a product of elements of the form $[v,h_{1}, \ldots, h_{m}]^{\varepsilon}$, with $m \geq 1$, $\varepsilon = \pm1$, $v \in {\cal{R}}$ and $h_{1}, \ldots, h_{m} \in F_{n+1}$. Since each $[R_{j}, F_{n+1}]$ is normal in $F_{n+1}$, we obtain the required claim. Hence $M = R [R_{2}, F_{n+1}] [R_{3}, F_{n+1} [R_{4}, F_{n+1}]$. 

\item Clearly $v \gamma_{\kappa + 3}(F_{n+1}) = [r, z_{1}, \ldots, z_{\kappa}]^{\varepsilon}\gamma_{\kappa + 3}(F_{n+1})$  with $\varepsilon = \pm 1$ and $z_{1}, \ldots, z_{\kappa} \in {\cal F}_{n+1}$. 
Without loss of generality,  we may assume $v = [r, z_{1}, \ldots, z_{\kappa}]$, with $r \in {\cal R}_{2}$ and $z_{1}, \ldots, z_{\kappa} \in {\cal F}_{n+1}$. 
In the next few lines, we consider the word $z_{1} \cdots z_{\kappa}$ as the $\kappa$-tuple $(z_{1}, \ldots, z_{\kappa})$. Let $1 \leq i_{1} < \cdots < i_{\nu} \leq \kappa$ be the $\nu$ different appearances of $t$ in the $\kappa$-tuple $(z_{1}, \ldots, z_{\kappa})$. 
By our hypothesis $\nu \geq 1$. 
Set $v(0) = [r,z_{1}, \ldots, z_{i_{1}-1}]$ with $z_{1}, \ldots, z_{i_{1}-1} \in {\cal F}_{n}$. For $i_{1} = 1$, $v(0) = r$. 
Write $g_{1} = r$, $g_{2} = z_{1}, \ldots,$ $g_{i_{1}} = z_{i_{1}-1}$. 
By Lemma \ref{le4.7}~(2) and working modulo $\gamma_{\kappa + 3}(F_{n+1})$  
$$
\begin{array}{rll}
v & = & -[[t, v(0)], z_{i_{1}+1}, \ldots, z_{\kappa}] \\
& = & \sum^{i_{1}}_{\lambda = 0}(\sum_{\sigma \in S_{i_{1},\lambda}}(-1)^{\lambda + 1}[t, g_{\sigma(1)}, \ldots, g_{\sigma(i_{1})}, z_{i_{1}+1}, \ldots, z_{\kappa}]).
\end{array}
$$
We point out that $r$ occurs in the word $g_{\sigma(1)}  \cdots g_{\sigma(i_{1})}$. By Lemma \ref{le4.7}~(1) (Jacobi identity), each of the above group commutator belongs to $[R_{3}, F_{n+1}, \ldots, F_{n+1}]$ with $\kappa - 1$ factors of $F_{n+1}$. By Proposition \ref{pro4.2}, $v \gamma_{\kappa + 3}(F_{n+1}) \in {\rm I}_{\kappa + 2}(M_{{\cal R}_{3}})$.  If there are no appearances of $t$  then,  by Proposition \ref{pro4.3}, $v\gamma_{\kappa + 3}(F_{n+1}) \in {\rm I}_{\kappa + 2}(M_{{\cal R}_{2}})$.                

\item We proceed as above. Without loss of generality, we may assume $v = [r, z_{1}, \ldots, z_{\kappa}]$, with $r \in {\cal R}_{4}$ and $z_{1}, \ldots, z_{\kappa} \in {\cal F}_{n+1}$. 
We separate two cases. First we assume that $z_{1}, \ldots, z_{\kappa} \in {\cal F}_{n}$. 
Let $r = [[t,f_{a}], [t,f_{b}]] \in {\cal{R}}_{4}$. 
Notice that $[[t,f_{a}],z]$ and $[[t,f_{b}],z] \in {\cal R}_{3}$ for all $z \in {\cal F}_{n}$. 
Using the Jacobi identity in the form $[x,y,z] = [x,z,y] + [x,[y,z]]$, Proposition \ref{pro3.6}~(4), and Proposition \ref{pro4.2}, we have $v\gamma_{\kappa + 5}(F_{n+1}) \in {\rm I}_{\kappa + 4}(M_{{\cal R}_{3}})$. 
Thus we assume that at least one of $z_{1}, \ldots, z_{\kappa}$ is $t$. 
Let $1 \leq i_{1} < \cdots < i_{\nu} \leq \kappa$ be the $\nu$ different appearances of $t$ in the $\kappa$-tuple $(z_{1}, \ldots, z_{\kappa})$. 
By our hypothesis $\nu \geq 1$. If $i_{1} \geq 2$, then $z_{1} \in {\cal F}_{n}$ and so, using the above argument, we get $v\gamma_{\kappa + 5}(F_{n+1}) \in {\rm I}_{\kappa + 4}(M_{{\cal R}_{3}})$. 
Thus we assume that  $z_{1} = t$. If $\kappa = 1$, then by Proposition \ref{pro3.6}~(3) and Proposition \ref{pro4.6}, we get $v\gamma_{6}(F_{n+1}) \in {\rm I}_{5}(M_{{\mathcal{R}}_{4}})$. 
Thus we assume that $\kappa \geq 2$. Working modulo ${\rm I}_{\kappa + 4}(M_{{\cal R}_{3}})$, using the Jacobi identity in the form $[x,y,z] = [x,z,y] + [x,[y,z]]$ and Proposition \ref{pro4.6}, we obtain $v\gamma_{\kappa + 5}(F_{n+1}) \in {\rm I}_{\kappa + 4}(M_{{\cal R}_{4}}) \oplus {\rm I}_{\kappa + 4}(M_{{\cal R}_{3}})$. 
Thus, in any case, $v\gamma_{\kappa + 5}(F_{n+1}) \in {\rm I}_{\kappa + 4}(M_{{\cal R}_{4}}) \oplus {\rm I}_{\kappa + 4}(M_{{\cal R}_{3}})$. \qed       
\end{enumerate}

\vskip .120 in

%Let $F_{m}$ be a free group of finite rank $m \geq 2$, freely generated by the set ${\cal{F}}_{m} = \{f_{1}, \ldots, f_{m}\}$. 
A group commutator $c_{j}$ of weight $w(c_{j})$ is defined as follows: The group commutators of weight one are the elements of ${\cal{F}}_{n}$, and if $i \neq j$ and $c_{i}$ and $c_{j}$ are group commutators of weights $w(c_{i})$ and $w(c_{j})$, respectively, then $[c_{i}, c_{j}]$ is a group commutator of weight $w(c_{i}) + w(c_{j})$. 

\begin{lemma}\label{le4.9} Let $F_{n}$ be a free group of rank $n \geq 2$, freely generated by the set ${\cal{F}}_{n} = \{f_{1}, \ldots, f_{n}\}$. 
\begin{enumerate}

\item Let $G$ be a subgroup of $F_{n}$ generated by a set of weighted group commutators $g_{1}, \ldots, g_{r}$ such that $2 \leq w(g_{1}) < w(g_{2}) < \cdots < w(g_{r})$, and choose a positive integer $\kappa$. Then every element of $G$ can be expressed in a form $b^{e_{1}}_{1} b^{e_{2}}_{2} \cdots b^{e_{\mu}}_{\mu}~{\rm mod}~\gamma_{\kappa+1}(G)$, where $e_{1}, \ldots, e_{\mu}$ are integers and $b_{1}, \ldots, b_{\mu}$ are group commutators in $g_{1}, \ldots, g_{r}$ such that $2 \leq w(b_{1}) < w(b_{2}) < \cdots < w(b_{\mu}) = \kappa$.

\item Let ${\cal{S}}$ be a subset of $F_{n}$ consisting of weighted group commutators $v_{1}, \ldots, v_{s}$ of $w(v_{1})= \cdots =  w(v_{s}) \geq 2$ and let $H$ be the subgroup of $F_{n}$ generated by ${\cal{S}}$. Then, for  $\kappa \geq 2$,  every element of $[H,F_{n}]$  can be expressed in a form $c^{d_{1}}_{1} \cdots c^{d_{\kappa-1}}_{\kappa-1}~{\rm mod}[H,\prescript{}{\kappa}{F_{n}}]$, $d_{1}, \ldots, d_{\kappa-1}$ are integers and each $c_{\nu}$ is a product of weighted group commutators of the form $[v, z_{1}, \ldots, z_{\nu}]$ with $v \in {\cal{S}}$ and $z_{1}, \ldots, z_{\nu} \in {\cal{F}}_{n}$.  

\end{enumerate} 

\end{lemma}

\pf \begin{enumerate}

\item It follows from a result of P. Hall \cite{phall} (see, also,  \cite[Chapter 3, Theorem 3.1]{cmz}). 

\item We induct on $\kappa$, and let $\kappa = 2$. Since $[H,F_{n}]$ is generated by the elements $[u,f]$ with $u \in H$ and $f \in F_{n} \setminus \{1\}$, it is enough to show that each $[u,f]$ is written in the required form. 
Using the well known group  commutator identities (see, for example, proof of Lemma \ref{le4.8}~(1)), working modulo $[H,F_{n},F_{n}]$ and since $[H, F_{n}, F_{n}]$ is normal in $F_{n}$, we obtain the desired result. Thus we assume that our claim is valid for some $\kappa \geq 2$, and let $g \in [H,F_{n}]$. By our inductive argument, $g$ is written as $c^{d_{1}}_{1} \cdots c^{d_{\kappa-1}}_{\kappa-1}w$ with $w \in [H,~_{\kappa}F_{n}]$ and each $c_{\mu}$ is a product of weighted group commutators of the form $[v, z_{1}, \ldots, z_{\mu}]$ with $v \in {\cal{S}}$ and $z_{1}, \ldots, z_{\mu} \in {\cal{F}}_{n}$. Working modulo $[H, ~_{\kappa + 1}F_{n}]$, we have $w$ is a product of weighted group commutators of the form $[v, z_{1}, \ldots, z_{\kappa}]$ with $v \in {\cal{S}}$ and $z_{1}, \ldots, z_{\kappa} \in {\cal{F}}_{n}$. Hence any element $g \in [H, F_{n}]$ can be expressed in the required form.~\qed 

\end{enumerate}

%%%%%%%%%%%%%%%%%%%%%%%%%%%%%%%%%%%%%%
\subsection{Main result}\label{subsec4.2}
%%%%%%%%%%%%%%%%%%%%%%%%%%%%%%%%%%%%%%

Let $M$ be the normal closure of ${\cal{R}}$ in $F_{n+1}$ and ${\cal L}(M) = \bigoplus _{d \geq 2}{\rm I}_{d}(M)$. Since ${\cal R}_{L} = {\mathcal{R}}_{2,L} \cup {\mathcal{R}}_{3,L} \cup {\mathcal{R}}_{4,L} \subset {\rm I}_{2}(M) \oplus {\rm I}_{3}(M) \oplus {\rm I}_{4}(M)$ and ${\cal L}(M)$ is an ideal of ${\rm gr}(F_{n+1})$, we have $J \subseteq {\cal L}(M)$.

\begin{proposition}\label{pro4.10}
${\cal L}(M) = J$. 
\end{proposition}

\pf Write ${\cal L} = {\mathcal{L}}(M_{{\mathcal{R}}_{2}}) + {\mathcal{L}}(M_{{\mathcal{R}}_{3}}) + {\mathcal{L}}(M_{{\mathcal{R}}_{4}})$. By Proposition \ref{pro4.3}, Proposition \ref{pro4.2}, Proposition \ref{pro4.6} and Proposition \ref{pro3.6}~(1), we have ${\cal L}$ is the additive direct sum of the Lie subalgebras ${\mathcal{L}}(M_{{\mathcal{R}}_{2}})$, ${\mathcal{L}}(M_{{\mathcal{R}}_{3}})$ and ${\mathcal{L}}(M_{{\mathcal{R}}_{4}})$ of ${\rm gr}(F_{n+1})$ and ${\cal L} = J$. 
Thus, for $d \geq 2$, $J \cap {\rm gr}_{d}(F_{n+1}) = {\rm I}_{d}(M_{{\mathcal{R}}_{2}}) \oplus {\rm I}_{d}(M_{{\mathcal{R}}_{3}}) \oplus {\rm I}_{d}(M_{{\mathcal{R}}_{4}})$. 

We claim that $J \cap {\rm gr}_{d}(F_{n+1}) = {\rm I}_{d}(M)$ for all $d \geq 2$. Since $J \cap {\rm gr}_{2}(F_{n+1}) = {\rm I}_{2}(M_{{\mathcal{R}}_{2}}) = {\rm I}_{2}(M)$, we may assume that $d \geq 3$. Let $u \in M \cap \gamma_{d}(F_{n+1})$ such that  $u \notin \gamma_{d+1}(F_{n+1})$. That is, $u\gamma_{d+1}(F_{n+1})$ is a non-trivial element of ${\rm I}_{d}(M) \leq {\rm gr}_{d}(F_{n+1})$. 
By Lemma \ref{le4.8}~(1), $u$ is written as $u = h h_{2} h_{3} h_{4}$ with $h \in R$ and $h_{j} \in [R_{j},F_{n+1}]$, $j \in \{2,3,4\}$. We point out that, for a positive integer $\kappa$, $[R_{j}, ~_{\kappa}F_{n+1}] \subseteq \gamma_{j+\kappa}(F_{n+1})$ with $j \in \{2,3,4\}$. 
Next, we work in the free nilpotent group $F_{n+1,d} = F_{n+1}/\gamma_{d+1}(F_{n+1})$. 
By Lemma \ref{le4.9} (for $G = R$ with a generating set ${\cal{R}}$, $H = R_{j}$ and ${\cal{S}} = {\cal{R}}_{j}$), and using the collection process for arbitrary group elements on $u$ (see \cite[Chapter 3, Theorem 3.1 and Theorem 3.5]{cmz}), the element $\overline{u} = u \gamma_{d+1}(F_{n+1})$ is written in $F_{n+1,d}$ as $\overline{w}_{2}\overline{w}_{3} \cdots \overline{w}_{d}$, where each $w_{i}$ is a product of weighted group commutators each of which has weight $i$ and in each weighted group commutator at least one element of ${\cal{R}}$ occurs. 
Notice that, for $i \in \{2, \ldots, d\}$, $w_{i} \in M \cap \gamma_{i}(F_{n+1})$. Since $u \in \gamma_{d}(F_{n+1}) \setminus \gamma_{d+1}(F_{n+1})$ and since the center of $F_{n+1,d}$ equals  $\gamma_{d}(F_{n+1,d}) = {\rm gr}_{d}(F_{n+1})$, we have $\overline{u} = \overline{w}_{d}$. 
Hence, working modulo $\gamma_{d+1}(F_{n+1})$, $u = \prod u^{\prime}_{\mu}$ is a product of weighted group commutators of weight $d$ in $F_{n+1}$ where in each $u^{\prime}_{\mu}$ at least one element of ${\cal R}$ occurs. Fix $\mu$ and let $u^{\prime}_{\mu}$ be a weighted group commutator of weight $d$ in $F_{n+1}$ of the form  $w(a_{1}, \ldots, a_{\lambda})$ where each $a_{i}$ is a weighted group commutator of weight $< d$. We separate the following cases. 
\begin{enumerate}
\item If an element of ${\cal R}_{3}$ occurs in at least one of $a_{1}, \ldots, a_{\lambda}$, then, by Proposition \ref{pro4.2}, $u^{\prime}_{\mu}\gamma_{d+1}(F_{n+1}) \in {\rm I}_{d}(M_{{\cal R}_{3}}) = L^{d}_{{\rm grad}}({\cal{E}} \wr \Omega^{*})$. 

\item If only elements of ${\cal{R}}_{2}$ occur in each  $a_{1}, \ldots, a_{\lambda}$, then, by Lemma \ref{le4.8}~(2), $u^{\prime}_{\mu}\gamma_{d+1}(F_{n+1})$ $\in {\rm I}_{d}(M_{{\cal{R}}_{2}}) \oplus {\rm I}_{d}(M_{{\cal{R}}_{3}})$ and so, by Proposition \ref{pro4.3} and Proposition \ref{pro4.2}, $u^{\prime}_{\mu}\gamma_{d+1}(F_{n+1})$ $\in J \cap {\rm gr}_{d}(F_{n+1})$. 

\item Finally, if an element of ${\cal{R}}_{4}$ occurs in at least one of $a_{1}, \ldots, a_{\lambda}$, then, by Lemma \ref{le4.8}~(3), we get $u^{\prime}_{\mu}\gamma_{d+1}(F_{n+1}) \in {\rm I}_{d}(M_{{\cal R}_{4}}) \oplus {\rm I}_{d}(M_{{\cal R}_{3}})$ and so, by Proposition \ref{pro4.2} and Proposition \ref{pro4.6}, $u^{\prime}_{\mu}\gamma_{d+1}(F_{n+1}) \in J \cap {\rm gr}_{d}(F_{n+1})$.  
\end{enumerate}
Therefore in any case $\overline{u} \in J \cap {\rm gr}_{d}(F_{n+1})$. Thus, for all $d \geq 2$, $J \cap {\rm gr}_{d}(F_{n+1}) = {\rm I}_{d}(M)$ and so ${\cal L}(M) = J$. \qed 

\begin{lemma}\label{le4.11}
$J$ is a direct summand of ${\rm gr}(F_{n+1})$.
\end{lemma}

\pf By Eq. $(2.8)$, ${\rm gr}(F_{n+1}) = {\mathcal{L}}_{\rm gr}(F_{n}) \oplus L(\Omega) \oplus L(\mathcal{E} \wr \Omega^{*})$. By Corollary \ref{co4.2} (for ${\rm gr}(F_{n}) = {\mathcal{L}}_{\rm gr}(F_{n})$ and $I = J \cap {\mathcal{L}}_{\rm gr}(F_{n})$), Lemma \ref{le4.4}, Lemma \ref{le4.5} and Proposition \ref{pro3.6}~(1), we obtain $J$ is a direct summand of ${\rm gr}(F_{n+1})$. \qed

\vskip .120 in

In the following result we give a presentation of the Lie algebra of the Formanek-Procesi group with basis a raag.

\begin{theorem}\label{thm4.12}
For $n \geq 2$, ${\rm gr}(F_{n+1}/M) \cong {\rm gr}(F_{n+1})/J$ as Lie algebras in a natural way.
\end{theorem}

\pf Throughout this proof, we write $G = F_{n+1}/M$. Since $M \subseteq F^{\prime}_{n+1}$, we have $G/G^{\prime} \cong F_{n+1}/F^{\prime}_{n+1}$ and ${\rm gr}(G)$ is generated as a Lie algebra by ${\cal G} = \{f_{i}G^{\prime}, f_{n+1}G^{\prime}=tG^{\prime}: i \in \n\}$. 
Since ${\rm gr}(F_{n+1})$ is free on $\overline{{\cal F}}_{n+1} = \{y_{1}, \ldots, y_{n}, \t\}$, the map $\zeta$ from ${\rm gr}(F_{n+1})$ into ${\rm gr}(G)$ satisfying the conditions $\zeta(y_{i}) = f_{i}G^{\prime}$, $i \in \n$, and $\zeta(\t) = tG^{\prime}$, extends uniquely to a Lie algebra homomorphism. 
Since ${\rm gr}(G)$ is generated as a Lie algebra by ${\cal G}$, we have $\zeta$ is onto. Hence ${\rm gr}(F_{n+1})/{\rm Ker}\zeta \cong {\rm gr}(G)$ as Lie algebras. For $d \geq 1$, 
$$
\begin{array}{rll}
{\rm gr}_{d}(G) & \cong & (\gamma_{d}(F_{n+1}) \gamma_{d+1}(F_{n+1})M)/(\gamma_{d+1}(F_{n+1})M) \\
& \cong & \gamma_{d}(F_{n+1})/(\gamma_{d}(F_{n+1}) \cap (\gamma_{d+1}(F_{n+1})M)).
\end{array}
$$
Since $\gamma_{d+1}(F_{n+1}) \subseteq \gamma_{d}(F_{n+1})$, we have, by the modular law, 
$$
\begin{array}{rll}
\gamma_{d}(F_{n+1})/(\gamma_{d}(F_{n+1}) \cap (\gamma_{d+1}(F_{n+1})M)) & = & \gamma_{d}(F_{n+1})/(\gamma_{d+1}(F_{n+1})(M \cap \gamma_{d}(F_{n+1})) \\
& \cong & {\rm gr}_{d}(F_{n+1})/{\rm I}_{d}(M)
\end{array}
$$ 
and so, for $d \geq 1$, ${\rm gr}_{d}(G) \cong {\rm gr}_{d}(F_{n+1})/{\rm I}_{d}(M)$ as $\mathbb{Z}$-modules in a natural way. Since $J$ is the ideal of ${\rm gr}(F_{n+1})$ generated by ${\mathcal{R}}_{L}$, $\zeta({\mathcal{R}}_{L}) = 0$ and ${\rm Ker}\zeta$ is ideal of ${\rm gr}(F_{n+1})$, we get  $J \subseteq {\rm Ker}\zeta$. Therefore $\zeta$ induces a Lie algebra epimorphism $\overline{\zeta}$ from ${\rm gr}(F_{n+1})/J$ onto ${\rm gr}(G)$. Since ${\rm gr}(F_{n+1}) = \bigoplus_{d \geq 1}{\rm gr}_{d}(F_{n+1})$ and $J = \bigoplus_{d \geq 2}(J \cap {\rm gr}_{d}(F_{n+1}))$, we have 
$$
{\rm gr}(F_{n+1})/J = \bigoplus_{d \geq 1}({\rm gr}_{d}(F_{n+1}) + J)/J.
$$ 
Thus, for $d \geq 1$, $\overline{\zeta}$ induces $\overline{\zeta}_{d}$, say, a $\mathbb{Z}$-linear mapping from $({\rm gr}_{d}(F_{n+1}) + J)/J$ onto ${\rm gr}_{d}(G)$. 
By Proposition \ref{pro4.10}, we have ${\rm I}_{d}(M) = J \cap {\rm gr}_{d}(F_{n+1})$ for all $d \geq 1$ and so we get 
$$
{\rm gr}_{d}(G)) \cong {\rm gr}_{d}(F_{n+1})/(J \cap {\rm gr}_{d}(F_{n+1}))
$$ 
as $\mathbb{Z}$-modules in a natural way. But, for $d \geq 1$,  
$$
\begin{array}{rll}
({\rm gr}_{d}(F_{n+1}) + J)/J & \cong & {\rm gr}_{d}(F_{n+1})/(J \cap {\rm gr}_{d}(F_{n+1})).
\end{array}
$$
By Lemma \ref{le4.11}, $({\rm gr}_{d}(F_{n+1}) + J)/J$ is a free $\mathbb{Z}$-module of finite rank and so ${\rm Ker}\overline{\zeta}_{d}$ is a free $\mathbb{Z}$-module of finite rank for all $d$. Since $
{\rm rank}({\rm gr}_{d}(G)) = {\rm rank}({\rm gr}_{d}(F_{n+1})/J \cap {\rm gr}_{d}(F_{n+1}))$,   
we have ${\rm Ker}\overline{\zeta}_{d} = \{1\}$ and so $\overline{\zeta}_{d}$ is an isomorphism for all $d$. Since $\overline{\zeta}$ is an epimorphism and each $\overline{\zeta}_{d}$ is isomorphism, we have $\overline{\zeta}$ is a Lie algebra isomorphism. Hence we obtain the required result. \qed

\bigskip

It follows from the above proof the following. 

\begin{corollary}\label{co4.13}
For $d \geq 1$, ${\rm gr}_{d}({\rm FP(H))}$ is a free abelian group of finite rank.
\end{corollary}

%%%%%%%%%%%%%%%%%%%%%%%%%%%%%%%%%%%%%
\section{Residual nilpotence of ${\rm FP(H)}$}\label{sect5}
%%%%%%%%%%%%%%%%%%%%%%%%%%%%%%%%%%%%%

A group $G$ is called residually nilpotent if to each $g \neq 1$ in $G$ there corresponds a normal subgroup $K_{g}$ such that $G/K_{g}$ is nilpotent and $g \notin K_{g}$. An equivalent assertion is that $G$ is residually nilpotent if $\bigcap_{c \geq 1} \gamma_{c}(G) = \{1\}$. A group $G$ is called a Magnus group if $G$ is residually nilpotent and each quotient group $\gamma_{c}(G)/\gamma_{c+1}(G)$ is torsion-free. 

Let $\theta$ be a partial commutative relation on $\{1,\ldots,n\}$ and  $A = \langle t, x_{1}, \ldots, x_{n}: [x_{\kappa}, x_{\ell}] = 1; \kappa > \ell; (\kappa, \ell) \in \theta\}$, $C = \langle y_{1}, \ldots, y_{n}: [y_{\kappa}, y_{\ell}] = 1; \kappa > \ell;(\kappa, \ell) \in \theta\}$.  
By its definition, $A$ is the free product of the infinite cyclic group $\<t\>$ generated by $t$ with the right-angled Artin group (raag) generated by $x_{1}, \ldots, x_{n}$. 
Since $\<x_1,\ldots, x_n\>$ is a raag, we have $A$ is also raag.  
Let ${\rm Aut}(A)$ be the automorphism group of $A$ with group operation $\psi_{1} \psi_{2} = \psi_{1} \circ \psi_{2}$. 
For any generator $y_{i}$ of $C$, let $\varphi_{y_{i}}: \{t, x_{1}, \ldots, x_{n}\} \rightarrow A$ be a map satisfying the conditions $\varphi_{y_{i}}(t) = tx_{i}$ and $\varphi_{y_{i}}(x_{j}) = x_{j}$ for all $j \in \{1, \ldots, n\}$. 
Clearly each $\varphi_{y_{i}}$ extends to a homomorphism of $A$. 
Since the set $\{tx_{i}, x_{1}, \ldots, x_{n}\}$ generates $A$, we have $\varphi_{y_{i}}$ is a group epimorphism of $A$. 
Since $A$ is residually nilpotent (see \cite{wade}), we get $A$ is residually finite and so $A$ is hopfian. Therefore each $\varphi_{y_{i}} \in {\rm Aut}(A)$. 

\begin{lemma}\label{le1}
With the previous notations, let $w = w(y_{1}, \ldots, y_{n}) = y^{\mu_{1}}_{i_{1}} \cdots y^{\mu_{m}}_{i_{m}}$ be a reduced word in $C$, with $i_{1}, \ldots, i_{m} \in \{1, \ldots, n\}$ and $\mu_{1}, \ldots, \mu_{m}$ non-zero integers, and let $\varphi_{w} = \varphi_{y^{\mu_{1}}_{i_{1}}} \cdots \varphi_{y^{\mu_{m}}_{i_{m}}}$.  Then $\varphi_{w}(t) = tx^{\mu_{1}}_{i_{1}} \cdots x^{\mu_{m}}_{i_{m}}$ and $\varphi_{w}(x_{j}) = x_{j}$ for all $j \in \{1, \ldots, n\}$.   
\end{lemma}

\pf Since $\varphi_{y_{i}}(x_{j}) = x_{j}$ for all $j \in \{1, \ldots, n\}$, we have $\varphi_{w}(x_{j}) = x_{j}$ for all $j \in \{1, \ldots, n\}$. 
Since $\varphi_{y^{\k}_{i}}(t) = tx^{\k}_{i}$ and $\varphi_{y^{\k}_{i}}(x_{j}) = x_{j}$ for all $i, j \in \{1, \ldots, n\}$ and $\k \in \mathbb{Z}$, we obtain the required result. \qed 

\vskip .120 in

The map $\varphi: C \rightarrow {\rm Aut}(A)$ defined by $\varphi(w(y_{1}, \ldots, y_{n})) = \varphi_{w(y_{1}, \ldots, y_{n})}$ for all $w(y_{1}, \ldots, y_{n}) \in C$ is well defined. Indeed, let $w_{1}(y_{1}, \ldots, y_{n}) = w_{2}(y_{1}, \ldots, y_{n})$. By Lemma \ref{le1}, 
$$
\varphi_{w_{1}(y_{1}, \ldots, y_{n})}(t) = tw_{1}(x_{1}, \ldots, x_{n}) = tw_{2}(x_{1}, \ldots, x_{n}) = \varphi_{w_{2}(y_{1}, \ldots, y_{n})}(t).
$$
Since $\varphi_{w_{1}(y_{1}, \ldots, y_{n})}(x_{j}) = \varphi_{w_{2}(y_{1}, \ldots, y_{n})}(x_{j}) = x_{j}$ for all $j \in \{1, \ldots, n\}$, we have $\varphi_{w_{1}(y_{1}, \ldots, y_{n})} = \varphi_{w_{2}(y_{1}, \ldots, y_{n})}$. It is clear enough that $\varphi$ is a homomorphism. Hence $C$ acts on $A$. Form the semidirect product $G = A \rtimes_{\varphi} C$ of $A$ and $C$; that is, the set $A \times C$ endowed with the group operation given by $(a_{1}, w_{1}) (a_{2}, w_{2}) = (a_{1} \varphi_{w_{1}}(a_{2}), w_{1}w_{2})$. Identify any element $w$ of $C$ with $(1,w)$ and any element $a$ of $A$ with $(a,1)$. Write $C_{\varphi} = \{(1,w): w \in C\}$ and $A_{\varphi} = \{(a,1): a \in A\}$. Clearly $C_{\varphi} \cong C$ and $A_{\varphi} \cong A$ as groups. Then $G = A_{\varphi} C_{\varphi}$, $A_{\varphi}$ is normal in $G$, $G/A_{\varphi} \cong C_{\varphi}$ and $A_{\varphi} \cap C_{\varphi} =  \{(1,1)\}$. We identify the group $C$ with $C_{\varphi}$ and view $C$ as a subgroup of $G$. Likewise, we identify $A$ with $A_{\varphi}$ and view it as a normal subgroup of $G$. With these identifications, the action of $C$ on $A$ becomes the restriction of the conjugation action in $G$; that is, $\varphi_{w}(a) = waw^{-1}$. Any element $g$ of $G$ can be written in a unique way as a product $g = aw$ for some $a \in A$ and $w \in C$. The group $G$ has a generating set $\{t, x_{1}, \ldots, x_{n}, y_{1}, \ldots, y_{n}\}$ and defining relations 
$$ 
[x_{\kappa}, x_{\ell}] = [y_{\kappa}, y_{\ell}] = 1:~\kappa > \ell;~(\kappa, \ell) \in \theta, 
$$
$$
y_{i}ty^{-1}_{i} = tx_{i},~~~~~~ y_{i}x_{j}y^{-1}_{i} = x_{j} ~{\rm for~all}~ i, j \in \{1, \ldots, n\}.
$$
$$
y_i^{-1}ty_i=tx_i^{-1}, ~~~~~~ y_i^{-1}x_jy_i=x_j~{\rm for~all}~ i, j \in \{1, \ldots, n\}.
$$

Guaschi and Pereiro \cite{gp} associate to a semidirect $G = A \rtimes_{\varphi} C$ a sequence of subgroups of $A$ which plays a central role to analyze the lower central series $\{\gamma_{d}(G)\}_{d \geq 1}$ of $G$ (see also \cite{suciu}). 
This sequence, $L = \{L_{d}\}_{d \geq 1}$, is inductively defined by setting $L_{1} = A$ and letting $L_{d+1}$ be the subgroup of $A$ generated by  $[L_{d}, A]$, $[A, \gamma_{d}(C)]$ and $[L_{d}, C]$. 
That is, $L_{d+1} = \langle [L_{d},A], [A, \gamma_{d}(C)], [L_{d}, C]\rangle$.   
The sequence $\{L_d\}_{d\ge 1}$ is a descending normal series of $A$.
We denote by $X$ the subgroup of $A$ generated by the set $\{x_{1}, \ldots, x_{n}\}$. 
For non-negative integers $r$ and $s$, with $r \geq s + 2$, we write $R_{r,s}(A,X) = [A, X, ~_{s}X,~_{r-s-2}A]$ and $R_{r,0}(A,X) = [A,X,~_{r-2}A]$. 
For the next few lines, let $Y_{1} = \cdots = Y_{s} = X$ and $Y_{s+1} = \cdots = Y_{r-2} = A$. For any permutation $\pi \in S_{r-2}$, we write $R_{r,s,\pi}(A,X) = [A,X,Y_{\pi(1)}, \ldots, Y_{\pi(r-2)}]$. 
Furthermore we write $\Gamma_{r,s}(A,X) = \langle \bigcup\limits_{\pi \in S_{r-2}}R_{r,s,\pi}(A,X) \rangle$. 
Notice that $\Gamma_{r,s}(A,X) \subset R_{r,0}(A,X) =  \gamma_r(A)$ for all $r \geq s + 2$.

\begin{lemma}\label{le2}
Let $G= A\rtimes C$  as above.
\begin{enumerate}

\item For $r \geq 1$, $\langle \gamma_{r}(A), [A,\gamma_{r}(C)] \rangle = \gamma_{r}(A)$. 

\item For $d \geq 4$, $
\langle \gamma_{d}(A), [\gamma_{d-1}(A), C] \rangle \subseteq \langle \gamma_{d}(A), \Gamma_{d-1,1}(A,X) \rangle$.

\end{enumerate} 
\end{lemma}

\pf 
\begin{enumerate}

\item We show that $[A,\g_r(C)] \subseteq \g_r(A)$. We work modulo $\g_r(A)$.  Since $A$ is a free product, every element of $A$ has the form $t^{k_1}w_1t^{k_2}w_2\ldots t^{k_m}w_m$ where $k_i\in\Z$ and $w_i$ are words in $\<x_1,\ldots,x_n\>$. Since $\gamma_{r}(C)$ is the normal closure of the set of  commutators $[y_{i_{1}}, \ldots, y_{i_{r}}]$ in $C$,  it suffices to show the result for the generators of $\g_r(C)$. Let $[y_{i_1},\ldots,y_{i_r}] \in \g_r(C)$ and $g\in C$.  Then 
$$
[t^{k_1}w_1t^{k_2}w_2\cdots t^{k_m}w_m, g[y_{i_1},\ldots,y_{i_r}]g^{-1}]=
$$
$$
w_m^{-1}t^{-k_m}\cdots w_1^{-1}t^{-k_1}\cdot g[y_{i_1},\ldots,y_{i_r}]^{-1}g^{-1}\cdot t^{k_1}w_1\ldots t^{k_m}w_m\cdot g[y_{i_1},\ldots,y_{i_r}]g^{-1}.
$$
Using Lemma \ref{le1} we get
$$w_m^{-1}t^{-k_m}\ldots w_1^{-1}t^{-k_1}\cdot (t\bar{g}[x_{i_1},\ldots,x_{i_r}]^{-1}(\bar{g})^{-1})^{k_1}w_1\ldots (t\bar{g}[x_{i_1},\ldots,x_{i_r}](\bar{g})^{-1})^{k_m}w_m$$
where $\bar{g}$ are words in $A$.  Working modulo $\g_r(A)$, the above is trivial and so we obtain the required result. 

\item Since $\gamma_{d-1}(A)$ is the normal closure of the set of commutators of the form $[a,b, z_{1},$ $\ldots, z_{d-3}]$ in $A$ with $a, b, z_{1}, \ldots, z_{d-3} \in \{t, x_{1}, \ldots, x_{n}\}$, it suffices to prove our claim for the generators of $\g_{d-1}(A)$. Let $u$ be a generator of $\g_{d-1}(A)$.  So $u$ has one of the following forms $u=g[t,x_j,z_1,\ldots z_{d-3}]g^{-1}$ or $u=g[x_k,x_j,z_1,\ldots, z_{d-3}]g^{-1}$
with $g\in A$ and $z_{1}, \ldots, z_{d-3}\in\{ t,x_1,\ldots, x_n\}$. Let $w=y_{i_1}\ldots y_{i_m}\in C$.  Then
$$
[u,w]=u^{-1}w^{-1}uw=u^{-1}\cdot\bar{g} [t(x_{i_1}\ldots x_{i_m})^{-1},x_j,\bar{z_1},\ldots,\bar{z}_{d-3}](\bar{g})^{-1}
$$ 
or 
$$
[u,w]=u^{-1}w^{-1}uw=u^{-1}\cdot\bar{g} [x_k,x_j,\bar{z_1},\ldots,\bar{z}_{d-3}](\bar{g})^{-1}
$$ 
where $\bar{z_r}=z_r$ if $z_r\in\{x_1,\ldots, x_n\}$ or $\bar{z_r}=t(x_{i_1}\ldots x_{i_m})^{-1}$ if $z_r=t$ and $\overline{g}$ are words in $A$.  Recall the following group commutator identities  
$$
[xy,z] = [x,z]^{y}[y,z],\eqno(5.1)
$$
$$
[x,yz] = [x,z][x,y]^{z}, \eqno(5.2)
$$
$$
[x^{-1},y] = ([x,y]^{-1})^{x^{-1}}, \eqno(5.3)
$$
$$
[x,y^{-1}] = ([x,y]^{-1})^{y^{-1}}. \eqno(5.4)
$$
Suppose that at least one of the $z_{1}, \ldots, z_{d-3} \in \{x_1,\ldots,x_n\}$. Using the commutator identities $(5.3)$ and $(5.1)$, working modulo $\gamma_{d}(A)$ and since $\gamma_{d}(A)$ is normal in $A$, 
we have each $[a,b, z_{1}, \ldots, z_{d-3}]$, with $a, b, z_{1}, \ldots, z_{d-3} \in \{t, x_{1}, \ldots, x_{n}\}$, belongs to $\Gamma_{d-1,1}(A,X)$. Let $u = [t,x_{j},z_{1}, \ldots,z_{d-3}]$. 
Using the identities $(5.1)-(5.4)$, working modulo $\gamma_{d}(A)$ and $\gamma_{d}(A)$ is normal in $A$, we get 
$$
[u,w] = u^{-1} u \prod (\bar{g}_{r}[x_{\mu},x_{j},\bar{z}_{1}, \ldots, \bar{z}_{d-3}](\bar{g}_{r})^{-1})^{\varepsilon} = \prod (\bar{g}_{r}[x_{\mu},x_{j},\bar{z}_{1}, \ldots, \bar{z}_{d-3}](\bar{g}_{r})^{-1})^{\varepsilon}
$$
where $\mu \in \{i_{1}, \ldots, i_{m}\}$, $\varepsilon = \pm 1$ and $\bar{g}_{r}$ are words in $A$, and so 
$$
[u,w] \in \langle  \gamma_{d}(A), \Gamma_{d-1,1}(A,X)\rangle.
$$
Similar arguments may be applied to $u = [x_{k},x_{j},z_{1}, \ldots,z_{d-3}]$. Assume now that $z_r=t$ for $r=1,\ldots,d-3$.  Let $u=[t,x_{j},t,\ldots,t]$. Then  
$$
[u,w]=u\cdot \bar{g}[t(x_{i_{1}} \cdots x_{i_{m}})^{-1},x_{j},t(x_{i_1} \ldots x_{i_m})^{-1},\ldots, t(x_{i_1}\ldots x_{i_m})^{-1}](\bar{g})^{-1}$$
Using the identities $(5.1)-(5.4)$, working modulo $\g_d(A)$ and since $\gamma_{d}(A)$ is normal in $A$, we have that this last element 
is a product of the form 
$$
u^{-1}u\cdot \prod (\overline{g}_r[a,b,\bar{z_1},\ldots,\bar{z}_{d-3}]\overline{g}_r^{-1})^{\varepsilon} = \prod (\overline{g}_r[a,b,\bar{z_1},\ldots,\bar{z}_{d-3}]\overline{g}_r^{-1})^{\varepsilon}
$$
where $a, b \in \{t, x_{1}, \ldots, x_{n}\}$, $\varepsilon = \pm 1$ and  at least one of the $z_1, \ldots, z_{d-3}\in\{x_1,\ldots,x_n\}$. Hence $[u,w] \in \langle \gamma_{d}(A), \G_{d-1,1}(A,X)\rangle$. Similar arguments may be applied to $u = [x_{k},x_{j},$ $t, \ldots,t]$. Therefore we obtain the desired result.~\qed

\end{enumerate}

\begin{lemma}\label{le4}
With the previous notations  
\begin{enumerate}

\item $L_{2} = \langle \gamma_{2}(A), X\rangle$.

\item $L_{3} = \gamma_{2}(A)$. 

\item $L_{4} = \langle \gamma_{3}(A), \gamma_{2}(X) \rangle$.  
 
\end{enumerate}   
\end{lemma}

\pf \begin{enumerate}

\item Let $H_{2} = \langle \gamma_{2}(A), X \rangle$. Clearly $H_{2}$ is normal in $A$. 
By the definition, $L_{2} = \langle \gamma_{2}(A), [A,C] \rangle$. Since $[t, y^{-1}_{j}] = t^{-1} (y_{j}ty^{-1}_{j}) = t^{-1}(tx_{j}) = x_{j}$ for all $j \in \{1, \ldots, n\}$, we get $x_{1},  \ldots, x_{n} \in L_{2}$. Since $L_{2}$ is normal in $A$, we obtain $H_{2} \subseteq L_{2}$. Let $a \in A$. Working modulo $\gamma_{2}(A)$, we write $a = t^{m}x^{m_{1}}_{1} \cdots x^{m_{n}}_{n}a^{\prime}$ where $m, m_{1}, \ldots, m_{n}$ are integers and $a^{\prime} \in \gamma_{2}(A)$. Write $b = x^{m_{1}}_{1} \cdots x^{m_{n}}_{n}$. By using the identity $(5.1)$, we get, for $w \in C$, 
$$
[a,w] = [t^{m}b,w]^{a^{\prime}} [a^{\prime},w] = ([t^{m},w]^{b}[b,w])^{a^{\prime}}[a^{\prime},w].
$$
Since $[b,w] = 1$ and $A$ is normal in $G$, we have $[a,w] = [t^{m},w]^{b}a_{1}$ with $a_{1} \in \gamma_{2}(A)$. Working modulo $\gamma_{2}(A)$, $[a,w] = [t,w]^{m}a_{2}$ with $a_{2} \in \gamma_{2}(A)$. Since $[t,y_{j}] = x^{-1}_{j}$ for all $j \in \{1, \ldots, n\}$, we obtain $[a,w] = w^{-m}(x_{1}, \ldots, x_{n}) a_{3}$ with $a_{3} \in \gamma_{2}(A)$. Therefore $L_{2} \subseteq H_{2}$ and so $L_{2} = H_{2} = \langle \gamma_{2}(A), X \rangle$.        

\item By the definition, $L_{3} = \langle [L_{2},A], [A, \gamma_{2}(C)], [L_{2},C]\rangle$. By Lemma \ref{le4}~(1) and the action of $C$ on $X$ we have
$$
L_{3} = \langle \gamma_{3}(A), [X,A], [A, \gamma_{2}(C)], [\gamma_{2}(A),C] \rangle.
$$
Since $\gamma_{2}(A)$ is the normal closure of the set $\{[t,x_{k}], [x_{i}, x_{j}]: k, i, j \in \{1, \ldots, n\}, i\neq j\}$ in $A$, we have $[X,A] = \gamma_{2}(A)$. Thus $
L_{3} = \langle \gamma_{2}(A), [A, \gamma_{2}(C)], [\gamma_{2}(A),C] \rangle$. 
By Lemma \ref{le2}~(1), we have $
L_{3} = \langle \gamma_{2}(A), [\gamma_{2}(A),C] \rangle$. 
Using the identity $(5.1)$ and working modulo $\gamma_{3}(A)$, we have $L_{3} \subseteq \gamma_{2}(A)$.  Since $\gamma_{2}(A) \subseteq L_{3}$, we obtain $L_{3} = \gamma_{2}(A)$.

\item By the definition and by Lemma \ref{le4}~(2), $
L_{4} = \langle \gamma_{3}(A), [A,\gamma_{3}(C)], [\gamma_{2}(A),C]\rangle$. 
By Lemma \ref{le2}~(1) (for $r = 3$), $
L_{4} = \langle \gamma_{3}(A), [\gamma_{2}(A),C]\rangle$. Since $[[t,v],y^{-1}] = [x,v] v^{\prime}$ with $v^{\prime} \in \gamma_{3}(A)$, we have $[x,v] \in L_{4}$ for all $x, v \in X$. Since $L_{4}$ is normal in $A$ and $\gamma_{2}(X)$ is the normal closure of the set $\{[x_{i}, x_{j}]: i, j \in \{1, \ldots, n\}\}$ in $X$, we obtain the required result. \qed

\end{enumerate}

\vskip .120 in

Let $d \geq 4$. If $d$ is even, we write 
$$
P_{d} = \left\langle \gamma_{d}(A), \Gamma_{d-i,i}(A,X): i \in \{1, \ldots, \frac{d}{2}-1\}\right\rangle=
$$
$$
 \left\< \g_d(A), \G_{d-1,1}(A,X),\G_{d-2,2}(A,X),\ldots, \G_{\frac d2+1,\frac d2-1}(A,X)\right\>
 $$
and if $d$ is odd, we write 
$$
P_{d} = \left\langle \gamma_{d}(A), \Gamma_{d-i,i}(A,X), \gamma_{\frac{d-1}{2}+1}(X): i \in \{1, \ldots, \frac{d-1}{2}-1\}\right\rangle=
$$
$$
\left\<\g_d(A),\G_{d-1,1}(A,X),\G_{d-2,2}(A,X),\ldots, \G_{\frac{d-1}{2}+2,\frac{d-1}{2}-1}(A,X),\g_{\frac{d-1}{2}+1}(X)\right\>
$$
It is obvious from the definition that $P_d\le \g_{[\frac d2]+1}(A)$ for every $d \geq 4$, where $[\frac d2]$ is the integer part of $\frac d2$. Notice that if $d$ is odd, then $[\frac{d}{2}] = \frac{d-1}{2}$.

\vskip .120 in

\begin{proposition}\label{pro1a}
With the previous notations
\begin{enumerate}

\item For $d \geq 4$, $P_d$ is a normal subgroup of $A$.  

\item For $d \geq 5$, $L_{d} \subseteq P_{d-1}$.

\end{enumerate} 

\end{proposition}

\pf \begin{enumerate}

\item Let $d \geq 4$. To prove that $P_d$ is a normal subgroup of $A$, it is enough to show that $[P_d,A] \subseteq P_d$. We separate two cases. 

\begin{enumerate}

\item If $d$ is even, then, for $i=2,\ldots, \frac d2+1$, we obviously have 
$$
[\G_{d-i,i}(X,A),A]\subseteq \G_{d-i+1,i}(X,A)\subseteq \G_{d-(i-1),i-1}(X,A)\subseteq P_d.
$$
For $i=1$ 
$$
[\G_{d-1,1}(X,A),A]\subseteq \g_d(A)\subseteq P_d.
$$  

\item If $d$ is odd,  then again
$$
[\G_{d-i,i}(X,A),A]\subseteq \G_{d-i+1,i}(X,A)\subseteq \G_{d-(i-1),i-1}(X,A)\subseteq P_d
$$
for every $i=2,\ldots, \frac{d-1}{2}-1$.
For $i=1$,  
$$
[\G_{d-1,1}(A,X),A]\subseteq \g_{d}(A)\subset P_d.
$$
Finally, we show that $[\g_{\frac{d-1}{2}+1}(X),A] \subseteq P_{d}$.  Since $[A,X,\prescript{}{\frac{d-1}{2}}X,A]$ is a normal subgroup of $A$ we have $[\g_{\frac{d-1}{2}+1}(X),A]$ is contained in $[A,X,\prescript{}{\frac{d-1}{2}-1}X,A]$.
By the definition of $\Gamma_{\frac{d-1}{2}+2,\frac{d-1}{2}-1}(A,X)$, we have $
[A,X,~_{\frac{d-1}{2}-1}X,A] \subseteq \Gamma_{\frac{d-1}{2}+2,\frac{d-1}{2}-1}(A,X) \subseteq P_{d}$. 
%Since $X \subseteq A$, we get $[A,X,~_{\frac{d-1}{2}}X] \subseteq [A,X,~_{\frac{d-1}{2}-1}X,A]$. Clearly $[A,X,~_{\frac{d-1}{2}}X,A]$ is a normal subgroup of $A$. Working modulo $[A,X,~_{\frac{d-1}{2}}X,A]$ and using Lemma \ref{le4.7}~(2) (for $G = A$ and $Y = X$), we get 
%$$
%[A, \gamma_{\frac{d-1}{2}+1}(X)] \equiv [A,X,~_{\frac{d-1}{2}}X].
%$$
%Since $[A,X,~_{\frac{d-1}{2}}X] \subseteq \Gamma_{\frac{d-1}{2} + 2, \frac{d-1}{2} -1}(A,X)$, we have $[A, \gamma_{\frac{d-1}{2}+1}(X)] \subseteq P_{d}$.   
Hence, we have $[\g_{\frac{d-1}{2}+1}(X),A]\subseteq P_d$.

\end{enumerate}

\item We induct on $d$ and let $d = 5$. By the definition, $L_{5} = \langle [L_{4},A], [A, \gamma_{4}(C)], [L_{4},C]\rangle$. By Lemma \ref{le4}~(3),  $
L_{5} = \langle \gamma_{4}(A), [A, \gamma_{2}(X)], [A, \gamma_{4}(C)], [\gamma_{3}(A), C], [\gamma_{2}(X),C] \rangle$. 
By Lemma \ref{le1},  $[\g_2(X),C]=\{1\}$ and so $L_{5} = \langle \gamma_{4}(A), [A, \gamma_{2}(X)], [A, \gamma_{4}(C)], [\gamma_{3}(A),C] \rangle$. By Lemma \ref{le2} (for $r = d = 4$), $L_{5} \subseteq \langle \gamma_{4}(A), \Gamma_{3,1}(A,X), [A, \gamma_{2}(X)] \rangle$. 
Using the Jacobi identity and working modulo $\gamma_{4}(A)$, we have $[A,\g_2(X)]\subseteq P_{4}$ and so $L_{5} \subseteq P_{4}$. Suppose that our claim is valid for some $d \geq 5$ and let $d$ be even, that is, $d \geq 6$. Then, by the inductive hypothesis,  $L_{d} \subseteq P_{d-1}$.  By its definition, 
$L_{d+1} = \langle [L_{d},A], [A,\gamma_{d}(C)], [L_{d}, C]\rangle$ and so $$
L_{d+1} \subseteq \langle [P_{d-1},A], [A,\gamma_{d}(C)], [P_{d-1},C] \rangle.
$$
Since $P_{d-1} = \langle \gamma_{d-1}(A), \Gamma_{(d-1)-i,i}(A,X), \gamma_{\frac{d}{2}}(X): i \in \{1, \ldots, \frac{d}{2}-2\}\rangle$, we get
$$
\begin{array}{rll}
L_{d+1} & \subseteq & \langle \gamma_{d}(A), [\Gamma_{(d-1)-i,i}(A,X),A], [\gamma_{\frac{d}{2}}(X),A], [\gamma_{d}(C),A], \\  
& & [\gamma_{d-1}(A),C], [\Gamma_{(d-1)-i,i}(A,X),C], [\gamma_{\frac{d}{2}}(X),C]:         i \in \{1, \ldots, \frac{d}{2}-2\}\rangle.
\end{array}
$$
By Lemma \ref{le2} and the action of $C$ on $A$ (Lemma \ref{le1}) 
%and since $[\gamma_{\frac{d}{2}}(X),A] = [A, \gamma_{\frac{d}{2}}(X)]$, 
$$
\begin{array}{rll}
L_{d+1} & \subseteq & \widetilde{P}_{d} = \langle \gamma_{d}(A), \Gamma_{d-1,1}(A,X), [\Gamma_{(d-1)-i,i}(A,X),A], [\gamma_{\frac{d}{2}}(X),A], \\  
& & [\Gamma_{(d-1)-i,i}(A,X),C]: i \in \{1, \ldots, \frac{d}{2}-2\}\rangle.
\end{array}
$$
Recall that $P_{d} = \langle \gamma_{d}(A), \Gamma_{d-j,j}(A,X): j \in \{1, \ldots, \frac{d}{2}-1\}\rangle$ and claim that $\widetilde{P}_{d} \subseteq P_{d}$. 
%Notice that $\rho + 1 = d$ and $\rho - [\frac{\rho}{2}] = [\frac{\rho}{2}]+1 = \frac{d}{2}$. 
Obviously,  $\g_d(A), \G_{d-1,1}(X,A)\subseteq P_d$ and $[\Gamma_{(d-1)-i,i}(A,X),A] \subseteq \Gamma_{d-i,i}(A,X)\subseteq P_d$  
for $i =1,\ldots, \frac{d}{2}-2$. 
Since $[A,X,\prescript{}{\frac d2 -2}X,A]$ is a normal subgroup of $A$ we have $[\g_{\frac d2}(X),A]\subseteq [A,X,\prescript{}{\frac d2 -2}X,A]$. Since $[A,X,\prescript{}{\frac d2 -2}X,A]\subseteq \Gamma_{\frac d2 +2,\frac d2 -2}(A,X)$ we get  $[\gamma_{\frac{d}{2}}(X),A] \subseteq P_{d}$.  
%Since $X \subseteq A$, we get $[A,X,~_{\frac{d}{2}-1}X] \subseteq [A,X,~_{\frac{d}{2}-2}X,A]$. Working modulo $[A,X,~_{\frac{d}{2}-1}X,A]$ and using Lemma \ref{le4.7}~(2) (for $G = A$ and $Y = X$), we get 
%$$
%[A, \gamma_{\frac{d}{2}}(X)] \equiv [A,X,~_{\frac{d}{2}-1}X].
%$$
%Since $[A,X,~_{\frac{d}{2}-1}X,A] \subseteq \Gamma_{\frac{d}{2} + 2, \frac{d}{2} -2}(A,X)$, we have $[A, \gamma_{\frac{d}{2}}(X)] \subseteq P_{d}$.  

Next we show that, for $i \in \{1, \ldots, \frac{d}{2}-2\}$, $[\Gamma_{(d-1)-i,i}(A,X),C] \subseteq P_{d}$. 
Let $i = \frac{d}{2}-2$. 
Then $
[\Gamma_{(d-1)-i,i}(A,X),C] = [\Gamma_{\frac{d}{2}+1, \frac{d}{2}-2}(A,X),C]$. 
By the definition 
$$
\Gamma_{\frac{d}{2}+1,\frac{d}{2}-2}(A,X) = \left\langle \bigcup_{\pi \in S_{\frac{d}{2}-1}}R_{\frac{d}{2}+1,\frac{d}{2}-2,\pi}(A,X)\right\rangle.
$$
Notice that $R_{\frac{d}{2}+1,\frac{d}{2}-2,{\rm Id}}(A,X) = [A,X,~_{\frac{d}{2}-2}X,A]$ with ${\rm Id}$ is the identity element of $S_{\frac{d}{2}-1}$ and, for $\pi \in S_{\frac{d}{2}-1} \setminus \{{\rm Id}\}$, $R_{\frac{d}{2}+1,\frac{d}{2}-2,\pi}(A,X) = [A,X,Y_{\pi(1)}, \ldots, Y_{\pi(\frac{d}{2}-1)}]$ with $Y_{1} = \cdots = Y_{\frac{d}{2}-2} = X$ and $Y_{\frac{d}{2}-1} = A$. 
Furthermore $[\Gamma_{\frac{d}{2}+1,\frac{d}{2}-2}(A,X),A] \subseteq \Gamma_{\frac{d}{2}+2,\frac{d}{2}-2}(A,X) \subseteq P_{d}$. 
Clearly $[\Gamma_{\frac{d}{2}+1,\frac{d}{2}-2}(A,X),A]$ is a normal subgroup of $A$. 
Let $[u,w]\in [\Gamma_{\frac{d}{2}+1, \frac{d}{2}-2}(A,X),C]$ with $u\in \Gamma_{\frac{d}{2}+1, \frac{d}{2}-2}(A,X)$ and $w\in C$.  
Then $[u,w]=u^{-1}w^{-1}uw$ and working modulo $[\Gamma_{\frac{d}{2}+1,\frac{d}{2}-2}(A,X),A]$ we may write $w^{-1}uw=u\cdot\prod_i w_i$ where $w_i\in \G_{\frac d2+1,\frac d2-1}(A,X)\subseteq P_d$. 
Therefore $[u,w]\in \G_{\frac d2+1,\frac d2-1}(A,X)$ and so $[\Gamma_{\frac{d}{2}+1, \frac{d}{2}-2}(A,X),C]\subseteq P_{d}$. 
Thus we may assume that $i < \frac{d}{2}-2$ and so $d \geq 8$.  
By the definition
$$
\Gamma_{(d-1)- i,i}(A,X) = \left\langle \bigcup_{\pi \in S_{(d-1) - 2}}R_{(d-1) - i,i, \pi}(A,X) \right\rangle = \left\langle \bigcup_{\pi \in S_{d-3}}R_{(d-1) - i,i, \pi}(A,X) \right\rangle.
$$
Since $d \geq 8$, we get $d - 3 \geq 5$. For $\pi \in S_{d-3}$, 
$$
R_{(d-1) - i,i, \pi}(A,X) = [A,X,Y_{\pi(1)}, \ldots,Y_{\pi(d-3)}]
$$ 
with $Y_{1} = \cdots = Y_{i} = X$ and $Y_{i+1} = \cdots = Y_{d-3-i} = A$. Furthermore 
$$
[\Gamma_{(d-1)-i,i}(A,X),A] \subseteq \Gamma_{d-i,i}(A,X)  \subseteq P_{d}.
$$
Now, working backwards (or using inverse induction), we may show that every element $[u,w]$ with $u\in \G_{(d-1)-i,i}(A,X)$ and $w\in C$ can be written as a product of elements in $\G_{d-(i+1),i+1}(A,X)$ and so it belongs to $P_d$.  Therefore $\widetilde{P}_{d} \subseteq P_{d}$. Hence $L_{d+1} \subseteq P_{d}$. Similar arguments may be applied if $d$ in our inductive hypothesis is odd. \qed

\end{enumerate}

\vskip .120 in

\begin{corollary}\label{co1}
For $\kappa \geq 2$, $L_{2 \kappa + 1}, L_{2 \kappa + 2} \subseteq \gamma_{\kappa + 1}(A)$.
\end{corollary}

\pf By Proposition \ref{pro1a}~(2), we have $L_{2 \kappa + 1} \subseteq P_{2 \kappa}$ and $P_{2 \kappa + 2} \subseteq P_{2 \kappa + 1}$. But 
\begin{center}
$P_{2 \kappa} = \left\langle \gamma_{2 \kappa}(A), \Gamma_{2 \kappa -i,i}(A,X): i \in \{1, \ldots, \kappa - 1\} \right\rangle$  ~~{\rm and} 

$P_{2 \kappa + 1} = \left\langle \gamma_{2 \kappa +1}(A), \Gamma_{2 \kappa + 1 -i,i}(A,X), \gamma_{\kappa + 1}(X): i \in \{1, \ldots, \kappa - 1\} \right\rangle.$
\end{center} 
For $i = \kappa - 1$, $\Gamma_{\kappa + 1, \kappa - 1}(A,X) = [A,X, ~_{\kappa - 1}X]$. Since $\Gamma_{\kappa + 1, \kappa - 1}(A,X) \subseteq \gamma_{\kappa + 1}(A)$ and $X \subseteq A$, we obtain the required result. \qed 

\vskip .120 in

\begin{corollary}\label{co2}
$\bigcap_{d \geq 1} L_{d} = \{1\}$.
\end{corollary}

\pf Since $A$ is residually nilpotent, we obtain from Lemma \ref{le4} and Corollary \ref{co1} the desired result. \qed

\vskip .120 in

\begin{theorem}\label{th5.7}
Let $G = A \rtimes_{\varphi} C$ be the split extension of groups  as above. Then $G$ is residually nilpotent. 
\end{theorem}

\pf Let $g \in G \setminus \{1\}$. We need to find a normal subgroup $N_{g}$ of $G$ such that $g \notin N_{g}$ and $G/N_{g}$ is nilpotent. First we assume that $g \notin A$. Since $G/A \cong C$ and $C$ is residually nilpotent ($C$ is raag), there exists a normal subgroup $M_{g}$ containing $A$ such that $g \notin M_{g}$ and $G/M_{g}$ is nilpotent. Thus we may assume that $g \in A$. By Corollary \ref{co2}, there exists a positive integer $d$ such that $g \in L_{d} \setminus L_{d+1}$. By a result of Guaschi and Pereiro \cite[Theorem 1.1]{gp}, $\gamma_{r}(G) = L_{r} \rtimes \gamma_{r}(C)$ for all $r \geq 1$. We claim that $g \notin \gamma_{d+1}(G)$. To get a contradiction, we assume that $g \in \gamma_{d+1}(G)$. Since $\gamma_{d+1}(G) = L_{d+1} \rtimes \gamma_{d+1}(C)$, we get $g \in L_{d+1}$ which is a contradiction. Since $G/\gamma_{d+1}(G)$ is nilpotent, we obtain the desired result. \qed 

\begin{corollary}\label{co5.8}
${\rm FP(H)}$ is residually nilpotent.
\end{corollary}

\pf As in the introduction,  ${\rm FP(H)}$ admits a presentation with generating set $b, a_{1,1},\ldots,a_{1,n}$, $a_{2,1},\ldots, a_{2,n}$ and defining relations
$$
ba_{1,k}a_{2,k}b^{-1}=a_{2,k}, \ \ k\in \n,
$$
$$
[a_{i,r},a_{i,s}]=1,~~ i=1,2;\ \ r > s;  (r,s)\in\theta~~\text{and}
$$
$$
[a_{1,k},a_{2,l}]=1, \ \ k,l\in \n
$$
with $\theta$ a partial commutation relation on $\{1,\ldots,n\}$.
By rearranging the above we get ${\rm FP(H)}=A\rtimes C$ with $A=\<b, a_{1,i}, i=1,\ldots,n\>$ and $C=\<a_{2,i}, i=1,\ldots, n\}$.  The action of $C$ on $A$ is described by the relations $a_{2,i}^{-1}ba_{2,i}
=ba_{1,i}$ for $i=1,\ldots, k$ and $a_{2,i}^{-1}a_{1,j}a_{2,i}=a_{1,j}$ for all $i,j\in\{1,\ldots, n\}$.
Hence, by Theorem \ref{th5.7}, it follows that ${\rm FP(H)}$ is residually nilpotent. \qed

\small
\bigskip

\noindent V. Metaftsis, Department of Mathematics, University of
the Aegean, Karlovassi, 832 00 Samos, Greece. {\it e-mail:}
vmet@aegean.gr

\bigskip

\noindent A.I. Papistas, Department of Mathematics, Aristotle
University of Thessaloniki, 541 24 Thessaloniki, Greece. {\it
e-mail:} apapist@math.auth.gr


\begin{thebibliography}{9}

\small

\bibitem{baht} Yu.A. Bahturin, Identical Relations in Lie Algebras, Nauka, Moscow, 1985 (in Russian).
English translation: VNU Science Press, Utrecht, 1987.

\bibitem{bour} N. Bourbaki, Lie Groups and Lie Algebras Part I, Hermann, Paris, 1987 (Chapters 1-3).

\bibitem{charney} R. Charney, An introduction to right angled Artin groups,  {\it Geom. Dedicata} {\bf 125} (2007), 141--158.

\bibitem{cmz} A.E. Clement, S. Majewicz and M. Zyman, The Theory of Nilpotent groups, Birkh$\ddot{a}$user, 2017.

\bibitem{dk} G. Duchamp and D. Krob, The free partially commutative Lie algebra: bases and ranks, {\it  Adv. Math.} {\bf 95} (1992), 92--126.

\bibitem{dk2} G. Duchamp and D. Krob, The lower central series
of the free partially commutative group, {\it Semigroup Forum} {\bf 45} (1992), 385--394.

\bibitem{fp} E. Formanek and C. Procesi, The automorphism group of a free group is not linear,  {\it J. Algebra} {\bf 149} (1992), 494--499.

\bibitem{gp} J. Guaschi and C.M. Pereiro, Lower central and derived series of semidirect products and applications to surface braid groups, {\it J. Pure and Applied Algebra}, {\bf 224} (2020), 1--39.

\bibitem{hall} M. Hall, A basis for free Lie rings and higher commutators in free groups, {\it Proc. Amer. Math. Soc.}, {\bf 1} (1950), 575--581.

\bibitem{phall} P. Hall, Some sufficient conditions for a group to be nilpotent, {\it Illinois J.  Math.}, {\bf 2} (1958), 787--801.

%\bibitem{hupBla} B. Huppert and N. Blackburn, Finite groups II, Springer-Verlag, Berlin, Heidelberg, New York, 1982. 

\bibitem{labu} J.P. Labute, The determination of the Lie algebra associated to the lower central series of a group, {\it Trans. Amer. Math. Soc.}, {\bf 288} (1985), no. 1, 51--57. 

\bibitem{laz} M. Lazard, Sur les groupes nilpotents et les anneaux de Lie, {\it Ann. Sci. Ecole Norm. Sup.}, {\bf 71} (3) (1954), 101--190.


\bibitem{mks} W. Magnus, A. Karrass and D. Solitar, Combinatorial Group Theory: presentations of groups in terms of generators and relations, Dover Publications, 1978.

\bibitem{sev} I. Sevaslidou, The associated Lie algebra of an HNN-extension of a group, {\it J. Gen. Lie Theory Appl.}, accepted.

\bibitem{suciu} A. Suciu, Lower central series and split extensions. arXiv:2105.14129.

\bibitem{wade} R. Wade, The lower central series of a right-angled Artin group,  {\it Enseign. Math.} {\bf 61} (2015), no. 3-4, 343--371.

\end{thebibliography}
\end{document}